\documentclass[a4paper]{amsart}

\usepackage[english]{babel}
\usepackage{amsmath,amsthm,amsfonts,amssymb}

\newenvironment{dimo}{\begin{proof}[Proof]}{\end{proof}}
\newtheorem{teo}{Theorem}[section]
\newtheorem{pro}[teo]{Proposition}
\newtheorem{lem}[teo]{Lemma}
\newtheorem{cor}[teo]{Corollary}
\newtheorem{defi}[teo]{Definition}
\newtheorem*{oss}{Remark}

\newtheorem*{nota}{Notation}

\numberwithin{equation}{section}

\newenvironment{sis}{\left\{\begin{aligned}}{\end{aligned}\right.}

\newcommand{\Q}{\mathbb{Q}}
\newcommand{\Z}{\mathbb{Z}}
\newcommand{\C}{\mathbb{C}}
\newcommand{\z}{\zeta}
\newcommand{\p}{\mathfrak{p}}

\begin{document}

\title{Ramification groups and Artin conductors of radical extensions of $\Q$}
\author{Filippo Viviani}
\address{Universita' degli studi di Roma Tor Vergata\\ Dipartimento di
matematica\\ via della ricerca scientifica 1\\ 00133 Roma, Italy}
\email{viviani@mat.uniroma2.it}

\begin{abstract}
We study the ramification properties of the extensions 
$\Q(\z_m,\sqrt[m]{a})/\Q$
under the hypothesis that $m$ is odd and if $p\mid m$ than either $p\nmid
v_p(a)$ or $p^{v_p(m)}\mid v_p(a)$ ($v_p(a)$ and $v_p(m)$ are the
exponents with which $p$ divides $a$ and $m$).
In particular we determine the higher ramification groups of the 
completed extensions and the Artin conductors of the characters of their 
Galois group.\\
As an application, we give formulas for the $p$-adique 
valuation of the discriminant of the studied global extensions 
with $m=p^r$.
\end{abstract}

\maketitle

\section{Introduction}
In this paper we study the ramification properties (ramification
groups and Artin conductors) of the normal radical extensions of
$\Q$, namely of the fields of the form $\Q(\z_m,\sqrt[m]{a})$ 
($\z_m$ a primitive $m$-th rooot of unity, $a\in \Z$),
under the hypothesis: (1) $m$ odd; (2) if $p\mid m$ then either
$p\nmid v_p(a)$ or $p^{v_p(m)}\mid v_p(a)$. While the first
hypothesis is assumed for simplicity (many strange phenomenas
appear when $2\mid m$ as the examples of the second section
show), the second hypothesis is a technical hypothesis that
unfortunately we weren't able to overcome (we will explain
in a moment why).

The interest in the radical extensions of the rationals is due to
the fact that they are the simplest and the more explicit normal
fields other than the abelian fields, so they are the ``first''
extensions not classified by the class field theory.
They have been studied under several point of
view: Westlund (\cite{wes}) and Komatsu (\cite{kom})
determined integral bases for $\Q(\sqrt[p]{a})$
and $\Q(\z_p,\sqrt[p]{a})$, respectively. Velez and Mann (see \cite{man-vel},
\cite{vel3}, \cite{vel5}) studied the factorization of primes in
$\Q(\sqrt[m]{a})$ and Jacobson and Velez (\cite{jacvel}) determined
in complete generality the Galois group of $\Q(\z_m,\sqrt[m]{a})$
(many complications arise when $2\mid m$, the case that we for
simplicity avoid). Eventually, the algebraic properties of the radical 
extensions have been studied by Acosta and Velez \cite{aco}, Mora and 
Velez \cite{mor} and Velez \cite{vel4}.

Our work is oriented in two new directions: the calculation of the
ramification groups and the Artin conductor of the characters of
the Galois group (for their definition and properties we refer
to the chapter IV and VI of the Serre's book \cite{ser}).
Let us now explain briefly what are the methods
that we used to obtain these results.

To calculate the ramification groups we first complete our
extensions with respect to $p$-adic valutation reducing in this
way to study the ramification groups of $p\neq 2$ in the local
extension $\Q_p(\z_{p^r},\sqrt[p^r]{a})$. Our original hypothesis on
$a$ splits in the two cases: (1) $p\nmid a$; (2) $p\mid\mid a$
(i.e. $p\mid a$ but $p^2\nmid a$, or $v_p(a)=1$).\\
Then the successive step is to calculate the ramification groups
of the extensions
$\Q_p(\z_{p},\sqrt[p^{i-1}]{a})<\Q_p(\z_{p},\sqrt[p^i]{a})$
and we succeed in this by finding an uniformizer (i.e. an element
of valuation $1$) and letting the generator of the
cyclic Galois group act on it (see Theorems 5.6 and 6.4).\\
The final step is a long and rather involved process of induction
which uses the knowledge of the ramification groups of the
cyclotomic fields (see \cite[Chapter IV, section 4]{ser})
as well as the functorial properties of the inferior and superior
ramification groups (see \cite[Chapter IV, section 3]{ser}).\\
Let us make at this point three remarks about these results and
the method used:\\
(1) the results obtained (Theorem 5.8 and 6.6) show that in the
non-abelian case the break-numbers in the superior ramification
groups are not necessarily integers (while the theorem of
Hasse-Arf (see \cite[Chapter V, section 7]{ser}) tells that this
always happens in the abelian case).\\
(2) the original hypothesis on the power of $p$ that divides $a$
is necessary because only in this case we are able to find an
uniformizer for the extension
$\Q_p(\z_{p},\sqrt[p^s]{a})<\Q_p(\z_{p},\sqrt[p^{s+1}]{a})$.
If one finds an uniformizer also for the other cases than the
same method will give the desired ramification groups in complete
generality.\\
(3) a much more simple and short proof could be obtained if one
finds directly an uniformizer for the whole extension
$\Q_p(\z_{p^r},\sqrt[p^r]{a})$. Unfortunately we were unable to
find this.

The other direction of our work is the calculation of the Artin
conductor of the characters of the Galois group
$G:=Gal(\Q(\z_m,\sqrt[m]{a})/\Q)$ (again we first reduce to the
case $m=p^r$).\\
So our first result is the explicit determination of the
characters of $G$ (Theo\-rem 3.7) and in order to do that we construct 
characters
in two ways:\\
(1) restricting characters of $(\Z/p^r\Z)^*$ under the projection
$G \twoheadrightarrow (\Z/p^r\Z)^*$.\\
(2) making Frobenius induction of characters of $\Z/p^r\Z$ under
the inclusion $\Z/p^r\Z \hookrightarrow G$.\\
After this we determine the local Artin conductor of a character $\chi$
by looking at its restriction at the ramification groups (Theorems 5.13 and 6.12).

Now let us summarize the contents of the various sections of the
paper.\\
 In the second section we recall some known results on
radical extensions and we prove, by applying a theorem of Schinzel
(see \cite{sch}), that if $k$ is a field that doesn't contain any
non trivial
$m$-root of unity then the polynomial $x^m-a$ remains irreducible
over $k(\z_m)$ if it is irreducible over $k$ (we shall apply this
for $k=\Q, \Q_p$).\\
In the third section we calculate the
characters of the group $Gal(\Q(\z_m,\sqrt[m]{a}))$ after having
decomposed it according to the prime powers that divide $m$.\\
In the fourth section we treat the case of tamely ramified primes.
In particular we show that if $p\mid a$ but $p\nmid m$ then $p$ is
tamely ramified and we calculate its ramification index (Theorem
4.3). Moreover in the case $p\mid m$ we show that the wild part of
the ramification is concentrated in the $p$-part
$\Q(\z_{p^r},\sqrt[p^r]{a})$ (Theorem 4.4).\\
The last two sections are devoted to study the ramification of $p$
in $\Q(\z_{p^r},\sqrt[p^r]{a})$ in the two cases $p\nmid a$ and
$p\mid a$.
In particular we compute the ramification groups of $p$ and the
$p$-local Artin conductor of the characters found in the third
section.
Then, by applying the conductor-discriminat formula, we calculate
the power of $p$ which divides the discriminant.

The referee pointed out to me that in the article\\
\textsc{H. Koch and E. de Shalit}, \textit{Metabelian local class field 
theory}. J. reine angew. Math. {\bf478} (1996), 85-106 \\
the authors studied the ramification groups of the maximal metabelian 
extension of a local field (of which the radical extensions considered 
here are particular examples) and asked about the compatibility between 
their and mine results. Actually to compare the results, one should 
compute the image of the metabelian norm map of our extensions and it's 
not at all clear to me how one can perform such a calculation (actually 
also in classical local class field theory, to determine the 
conductor of an abelian extension it's often simpler to calculate the 
ramification groups or the Artin conductors than the image of the norm 
map). The advantage of my results is that they are explicit and permit to 
avoid this much more general and elaborated theory.   
 
\emph{Aknowledgement}: This work is the result of my master thesis
which was made at the University of Pisa under the direction of
prof. Roberto Dvornicich which we thank for his supervision and
encouragements.

\section{Some results on radical extensions}

In this section we collect some known results on radical extensions that
we shall need in the next sections. We shall always consider the equation
$x^m-a$ defined over a field $k$ such that $\mbox{char}(k)\nmid m$ and we
shall restrict ourselves to the case in which $m$ is odd (when the prime
$2$ appears in the factorization of $m$, new strange phenomenas occur so
that for semplicity we prefer to avoid these complications).

\begin{teo}
The equation $x^m-a$ (with $m$ odd) is irreducible if and only if
$a\not\in k^{p}$ for every $p\mid m$.
\end{teo}
\begin{dimo}
See \cite[Chapter VI, Theorem 9.1]{lang}. 
\end{dimo}
\begin{teo}
Let $x^m-a$ irreducible over $k$ with $2\nmid m$. Then $k(\sqrt[m]{a})/k$
has the unique subfield property, i.e. for every divisor $d$ of $[L:K]$
there exists a unique intermediate field $M$ such that $[M:K]=d$.\\
Precisely, if $d\mid m$, then the unique subextension $L$ of degree $d$
over $k$ is $L=k(\sqrt[d]{a})$.
\end{teo}
\begin{dimo}
See \cite[Theorem 2.1]{aco}.
\end{dimo}
\begin{teo}
Let $x^m-a$ and $x^m-b$ irreducible over $k$, with $m$ odd. If
$k(\sqrt[m]{a})=k(\sqrt[m]{b})$, then $\sqrt[m]{b}=c(\sqrt[m]{a})^t$ for
some $c\in k$ and $t\in \mathbb{N}$ such that $(t,m)=1$.
\end{teo}
\begin{dimo}
It follows from the preceding theorem and from \cite[Lemma 2.3]{vel2}.
\end{dimo}
\begin{oss}
All these three results are false if $m$ is divisible for $2$ as the
following examples show:
\begin{enumerate}
\item $x^4-(-4)=(x^2+2x+2)(x^2-2x+2)$ but $-4\not\in \Q^2$;
\item $x^4+1$ is irreducible over $\Q$ but $\Q(\sqrt[4]{-1})=\Q(\zeta_8)
=\Q(i,\sqrt{2})$ has three subfields of degree $2$: $\Q(i)$,
$\Q(\sqrt{2})$, $\Q(\sqrt{-2})$;
\item
$\Q(\sqrt[8]{-1})=\Q(\sqrt{2},i,\sqrt{\sqrt{2}+2})=\Q(\sqrt[8]{-16})$
but $16\not\in\Q^8$.
\end{enumerate}
\end{oss}
The last result we need is a theorem of Schinzel characterizing the
abelian radical extensions, i.e. those radical extensions whose
normal closure has abelian Galois group.
 \begin{teo}[Schinzel]
Let $k$ be a field and let $m$ be a natural such that $\mbox{char}(k)\nmid
m$. Denote with $\omega_m$ the number of the $m$-roots of unity
contained in $k$.\\
Then the Galois group of $x^m-a$ over $k$ is abelian (i.e.
$\mbox{Gal}(k(\z_m,\sqrt[m]{a})/k)$ is abelian) if and only if
$a^{\omega_m}=\gamma^m$ for some $\gamma\in k$.
\end{teo}
\begin{dimo}
See the original paper of Schinzel (\cite[Theorem 2]{sch}). For other
proofs see \cite[Theorem 2.1]{vel2} or \cite[Lemma 7]{woj}, while a 
nice generalization of this theorem is contained in \cite{vel1}.
\end{dimo}
Using the theorem of Schinzel we can prove the following proposition.
\begin{pro}
Let $1\leq m\mid n$ and let $k$ be a field such that $\mbox{char}(k)\nmid
m$. If $k$ doesn't contain any $m$-root of unity other than the
identity then an element $a$ of $k$ is a $m$-power in
$k(\z_n)$ if and only if it is a $m$-power in $k$.
\end{pro}
\begin{dimo}
The ``if''  part is obvious. Conversely assume that $a\in k(\z_n)^m$. Then
$$k(\sqrt[m]{a})\subset k(\z_n)\Rightarrow k(\sqrt[m]{a},\z_m)\subset
k(\sqrt[m]{a},\z_n)\subset k(\z_n)$$
and so $x^m-a$ has abelian Galois group over $k$. But then the theorem
of Schinzel implies $a\in k^m$, q.e.d.
\end{dimo}
\begin{oss}
The preceding result is false if the field contains some $m$-root of
unity other than the identity as the following example shows:\\
$-1$ is not a square in $\Q$ but it becomes a square in
$\Q(\z_4)=\Q(i)$ (the reason is that $\Q$ contains the non trivial
$2$-root of unity $-1$).
\end{oss}
We can now apply these results on radical extensions to the situation we
are concerned with, i.e. the irreducibility of the polynomial $x^m-a$
defined over $\Q$ and with $m$ odd.
\begin{cor}
If $m$ is odd and $a\not\in\Q^p$ for every $p\mid m$, then the polynomial
$x^m-a$ is
irreducible over $\Q$ and so the extension $\Q(\sqrt[m]{a})/\Q$ has
degree $m$.
\end{cor}
\begin{dimo}
It follows at once from Theorem 2.1.
\end{dimo}
Moreover in the next sections we will consider the normal closure of
$\Q(\sqrt[m]{a})$, i.e $\Q(\z_m,\sqrt[m]{a})$. The next result tells us
what is the degree of this extension.
\begin{cor}
If $m$ is odd and $a\not\in\Q^p$ for every $p\mid m$, then $x^m-a$
is irreducible over $\Q(\z_m)$ and so
$[\Q(\z_m,\sqrt[m]{a}):\Q]=\phi(m)m$. \end{cor}
\begin{dimo}
It follows at once from Proposition 2.5 and Corollary 2.6 after
observing that $\Q$ doesn't contain any $m$-root of unity other than
$1$ if $m$ is odd.
\end{dimo}
\begin{oss}
Again the preceding result is false if $2\mid m$ as the "usual"
example shows:\\
$x^4+1$ is irreducible over $\Q$ but over $\Q(\z_4)=\Q(i)$ it splits as
$x^4+1=(x^2+i)(x^2-i)$ and hence
$[\Q(\sqrt[4]{-1},\z_4):\Q]=4<\phi(4)4=8$.\\
For an analysis of the degree of the splitting field of the polynomial
$x^{2^s}-a$ as well as of its Galois group see the paper
\cite{jacvel}.
\end{oss}
We end this preliminary section with this useful splitting result.
\begin{pro}
If $m=\prod_{i=1}^s p_i^{r_i}$ then the extension $k(\sqrt[m]{a})$ is the
compositum of the extensions $k(\sqrt[p_i^{r_i}]{a})$ for $i=1,\cdots,s$,
i.e.
$$k(\sqrt[m]{a})=k(\sqrt[p_1^{r_1}]{a})\cdots k(\sqrt[p_s^{r_s}]{a}).$$
\end{pro}
\begin{dimo}
It's enough to prove that if $m=m_1m_2$ with $m_1$ and
$m_2$ relatively prime, then
$k(\sqrt[m]{a})=k(\sqrt[m_1]{a})k(\sqrt[m_2]{a})$.
The inclusion $k(\sqrt[m]{a})\supset k(\sqrt[m_1]{a})k(\sqrt[m_2]{a})$
is obvious. On the other hand, since $(m_1,m_2)=1$, there exist $s,t\in\Z$ such that
$sm_1+tm_2=1$. But this imply $(\sqrt[m_1]{a})^t(\sqrt[m_2]{a})^s=
(\sqrt[m_1m_2]{a})^{tm_2+sm_1}=\sqrt[m]{a}$, q.e.d.
\end{dimo}

\section{Characters of the Galois groups of $x^m-a$}

First of all we want to describe the Galois group
$\mbox{Gal}(\Q(\z_m,\sqrt[m]{a})/\Q)$.
\begin{defi}
The holomorphic group of a finite group $G$ (non necessarily abelian, 
although we use the addittive notation) is the semidirect product
of $G$ with $\text{Aut}(G)$ (indicated with $G\rtimes\text{Aut}(G)$), that
is the set of pairs $\{\,(g,\sigma)$\,:\,$g\in G$,
 $\sigma\in\text{Aut}(G)\}$ with the multiplication given by
\begin{eqnarray}
(g,\sigma)(h,\tau)=(g+\sigma(h),\sigma\circ\tau). 
\end{eqnarray}
\end{defi}
\begin{nota}
We shall denote by $C(m)$ the cyclic group of order $m$
(identified with $\Z/m\Z$)
and with $G(m)$ the group of its automorphisms (identified with
$(\Z/m\Z)^*$). We shall denote the holomorphic group of $C(m)$
with $K(m):=C(m)\rtimes G(m)$ (the letter $K$ stands
for Kummer who first studied this kind of extensions).
\end{nota}
\begin{pro}
Suppose that $x^m-a$ is irreducible over $\Q$. Then the Galois group
of $x^m-a$ is isomorphic to the holomorphic group of the cyclic group of
order $m$, i.e.
$$\mbox{Gal}(\Q(\z_m,\sqrt[m]{a})/\Q)\cong C(m)\rtimes G(m)=K(m).$$
\end{pro}
\begin{dimo}
Every element $\sigma$ of $\mbox{Gal}(\Q(\z_m,\sqrt[m]{a})/\Q)$ is
uniquely determined by its values on the generators of the extension
and on these it must hold
\begin{alignat*}{2}
&\sigma(\sqrt[m]{a})=\z_m^i\sqrt[m]{a}&\qquad&i\in C(m),\\
&\sigma(\z_m)=\z_m^k&&k\in G(m).
\end{alignat*}
Then we can define an injective homomorphism
\begin{alignat*}{2}
\Gamma:\mbox{Gal}(\Q(\z_m,\sqrt[m]{a})/&\Q) \longrightarrow&
       &C(m)\rtimes G(m)\\
       &\sigma\stackrel{\Gamma}{\longmapsto}   & & (i,k).
\end{alignat*}
But this is an isomorphism since both groups
have cardinality equal to $m\phi(m)$ (see Proposition 2.7).
\end{dimo}
Next a useful splitting result.
\begin{pro}
If $m=\prod_i{p_i^{r_i}}$ then $K(m)\cong \prod_iK(p_i^{r_i})$.
\end{pro}
\begin{dimo}
This follows easily from the analogue property of the
groups $C(m)$ and $G(m)$.
\end{dimo}
So we have reduced ourselves to study the characters of the Kummer group
$K(p^r)=C(p^r)\rtimes G(p^r)$, for $p$ odd prime.
\begin{nota}
In what follows we shall adopt the following convention:
the elements of $K(p^r)$ will be denoted with $z^i\sigma$,
where the roman letters $i,j,k,\cdots$  will indicate elements of $C(p^r)$,
the greek letters $\sigma,\tau,\cdots$
will indicate elements of $G(p^r)$ and the letter $z$ is an
auxiliary letter that will allow to transform the multiplicative
notation for $K(p^r)$ into the additive notation for its subgroup
$C(p^r)$.
\end{nota}
With this notation, the product in $K(p^r)$ is ruled by the following equation
\begin{eqnarray}
z^i\sigma z^j\tau=z^{i+\sigma j}\sigma\tau. 
\end{eqnarray}
Observe that $G(p^r)$ acts on the normal subgroup $C(p^r)$ by
\begin{eqnarray}
\sigma z^i\sigma^{-1}=z^{\sigma i}. 
\end{eqnarray}
Let us determine the conjugacy classes of $K(p^r)$.
\begin{teo}
Let $z^i\sigma$ be an element of $K(p^r)$ ($p\neq 2$) and let
$$\begin{cases}
\alpha=v_p(\sigma-1)\qquad &0\leq\alpha\leq r,\\
\beta=v_p(i)        \qquad &0\leq \beta \leq r.\\
\end{cases}$$
Then the conjugacy class of $z^i\sigma$ is
$$[z^i\sigma]=\begin{sis}
               &\{\,z^j\sigma\,:\,v_p(j)=v_p(i)=\beta\,\}\quad&&\text{ if }
                               &&0\leq\beta<\alpha\\
               &\{\,z^j\sigma\,:\,v_p(j)\geq v_p(\sigma-1)=\alpha\,\}\quad
                            &&\text{ if }&&\alpha\leq\beta.
              \end{sis}$$
\end{teo}
\begin{dimo}
It's enough to consider the conjugates by elements of $C(p^r)$ and
$G(p^r)$
\begin{align*}
&z^{-k}z^i\sigma z^k=z^{i+k(\sigma-1)}\sigma \tag{*},\\
&\tau z^i\sigma\tau^{-1}=z^{i\tau}\sigma \tag{**}.\\
\end{align*}
\fbox{ $0\leq\beta<\alpha$ } Let us prove the two inclusions in the
                             statement of the theorem.

\underline{$\subset$} In (*) we have $v_p(i+k(\sigma-1))=v_p(i)$
                      since $v_p(i)<v_p(\sigma-1)\leq v_p(k(\sigma-1))$.
In (**) $v_p(i\sigma)=v_p(i)+v_p(\sigma)=v_p(i)$.

\underline{$\supset$} Let $j\in C(p^r)$ be such that $v_p(j)=v_p(i)$.
Then $\tau:=ji^{-1}$ has $p$-adic valuation equal to $0$ and so it belongs
to $G(p^r)$. Then from (**) we see that $z^j\sigma\in[z^i\sigma]$.\\

\fbox{ $\alpha\leq\beta$ } Let us prove the two inclusions.

\underline{$\subset$} In (*) we have $v_p(i+k(\sigma-1))\geq \min\{v_p(i),
\,v_p(k)+v_p(\sigma-1)\}\geq\min\{\alpha,\,v_p(k)+\alpha\}=\alpha$.
In (**) $v_p(i\tau)=v_p(i)=\beta\geq\alpha$.

\underline{$\supset$} Given $j\in C(p^r)$ such that  $v_p(j)\geq v_p(\sigma-1)$,
the equation $j=i+k(\sigma-1)$ is solvable for some $k\in C(p^r)$
and so from (*) we conclude $z^j\sigma\in [z^i\sigma]$.
\end{dimo}
Recall that the group $G(p^r)$ has a filtration given by the
subgroups $G(p^r)^{\alpha}=\{\,\sigma\in
G(p^r)\,:\,v_p(\sigma-1)\geq\alpha\,\}$.
\begin{cor}
Given $\sigma\in G(p^r)$ such that $v_p(\sigma-1)=\alpha$ (which from now
on we shall denote with $\sigma_{\alpha}$), the set
$\{\,z^i\sigma\,:\,i\in C(p^r)\}$
is invariant under conjugacy and splits in the $\alpha+1$ classes
$$\begin{aligned}
&[z\sigma_{\alpha}]&&=&&\{\,z^j\sigma_{\alpha}\,:\,v_p(j)=0\}\\
&\;\cdots   && &&\qquad\; \cdots                     \\
&[z^{p^i}\sigma_{\alpha}]&&=&&\{\,z^j\sigma_{\alpha}\,:\,v_p(j)=i\}\\
&\;\cdots   && &&\qquad\; \cdots                     \\
&[z^{p^{\alpha-1}}\sigma_{\alpha}]&&=&&\{\,z^j\sigma_{\alpha}\,:\,v_p(j)=\alpha-1\}\\
&[z^{p^{\alpha}}\sigma_{\alpha}]&&=&&\{\,z^j\sigma_{\alpha}\,:\,v_p(j)\geq\alpha\}.\\
\end{aligned}$$
\end{cor}
Now we can count the number of conjugacy classes of $K(p^r)$.
\begin{pro}
The number of conjugacy classes of $K(p^r)$ ($p$ odd prime) is equal to
$$\#\{\mbox{Conjugacy classes}\}=(p-1)p^{r-1}+\frac{p^r-1}{p-1}.$$
\end{pro}
\begin{dimo}
According to the preceding corollary we have
$$\begin{aligned}
&\#\{[z\sigma]\}=p^{r-1}(p-1)\\
&\#\{[z^p\sigma]\,:\,[z^p\sigma]\neq [z\sigma]\}=\#\{\sigma\,:\,
v_p(\sigma-1)\geq 1\}=|G(p^r)^1|=p^{r-1}\\
&\quad \cdots\\
&\#\{[z^{p^r}\sigma]\,:\,[z^{p^r}\sigma]\neq[z^{p^{r-1}}\sigma]\}=\#\{\sigma\,:\,
v_p(\sigma-1)\geq r\}=|G(p^r)^r|=1.\\
\end{aligned}$$
So the number of conjugacy classes is
$$\#\{\mbox{Conjugacy classes}\}=p^{r-1}(p-1)+p^{r-1}+\cdots+p^{r-r}=
p^{r-1}(p-1)+\frac{p^r-1}{p-1}.$$
\end{dimo}
Before we determine the characters of $K(p^r)$, we recall some
facts about the characters of the group $G(p^r)=(\Z/p^r\Z)^*$, with $p$
odd prime. First of all we know that
$$G(p^r)=<g^{p^{r-1}}>\oplus<1+p>\cong C(p-1)\times C(p^{r-1})$$
with $0<g<p$ a generator of the cyclic group $G(p)$.\\
Besides recall that $G(p^r)$ has a natural filtration
$$G(p^r)\supset G(p^r)^1\supset\cdots\supset G(p^r)^{r-1}\supset G(p^r)^r=
\{1\}$$
where $G(p^r)^k=\{\sigma\in G(p^r)\,:\,v_p(\sigma-1)\geq k\}=<1+p^k>\cong
C(p^{r-k})$, for $1\leq k\leq r$, and moreover we have
$G(p^r)/G(p^r)^k=G(p^k)$.\\
If we translate this information at the level of characters we obtain:

(1) $G(p^k)^*\subset G(p^r)^*$ through the
    projection $G(p^r)\twoheadrightarrow G(p^k)$;

(2) $G(p^r)^*/G(p^k)^*=(G(p^r)^k)^*\cong C(p^{r-k})$ through the
    inclusion $G(p^r)^k=$\\
    $<1+p^k>\hookrightarrow G(p^r)$.

\begin{nota}
In what follows we shall denote the characters of $G(p^r)$ with $\psi^r$,
and with $\psi^r_k$ a fixed system of representatives
for the lateral cosets of $G(p^k)^*$ in $G(p^r)^*$, in such a
way that, when restricted, they give all the characters of $G(p^r)^k$.
\end{nota}
With these notation we can now determine all the characters of $K(p^r)$.
\begin{teo}
The irreducible characters of $K(p^r)$ (for $p$ odd prime) are
\begin{center}
\begin{tabular}{||p{8.5cm}||c|c||}
\hline
\hline
\bfseries CHARACTERS       & Number       & Degree\\
\hline
\hline
$\psi^r\qquad\text{ with }\psi^r\in G(p^r)^*$       &$p^{r-1}(p-1)$&$1$\\
\hline
$\psi^r_r\otimes\chi_r^r\;\text{with }\psi^r_r\mbox{ syst. of repr. for }
G(p^r)^*/G(p^r)^*$
& $1$ & $p^{r-1}(p-1)$\\
\hline
$\cdots$                  &$\cdots$      &$\cdots$ \\
\hline
$\psi^r_k\otimes\chi^r_k\;\text{with }\psi^r_k\mbox{ syst. of repr. for }
G(p^r)^*/G(p^k)^*$
& $p^{r-k}$ & $p^{k-1}(p-1)$\\
\hline
$\cdots$                 &$\cdots$      &$\cdots$ \\
\hline
$\psi^r_1\otimes\chi^r_1\;\text{with }\psi^r_1\mbox{ syst. of repr. for }
G(p^r)^*/G(p^1)^*$
& $p^{r-1}$ &$(p-1)$\\
\hline
\hline
\end{tabular}
\end{center}
where

(i) $\otimes$ means the tensorial product of representations which, at
the level of characters, becomes pointwise product;

(ii) $\psi^r$ is the character defined by
    $$\psi^r(z^{p^{\beta}}\sigma_{\alpha})=\psi^r(\sigma_{\alpha})$$
    that is the character induced on $K(p^r)$ from $G(p^r)$ through the
    projection $K(p^r)\twoheadrightarrow G(p^r)$. Analogously
    the $\psi^r_k$ are seen as characters on $K(p^r)$.

(iii) $\chi^r_k$, $1\leq k\leq r$, is the character defined by
$$\chi^r_k([z^{p^{\beta}}\sigma_{\alpha}])=\left\{\begin{aligned}
&0\quad&&\text{ if }&&\alpha<k&&\text{ or }&&\beta<k-1,\\
&-p^{k-1}&&\text{ if }&&k\leq\alpha&&\text{ and }&&\beta=k-1,\\
&p^{k-1}(p-1)&&\text{ if }&&k\leq\alpha&&\text{ and }&&k-1<\beta.\\
\end{aligned}\right.$$\\
(Recall that $\sigma_{\alpha}$ indicates an element of $K(p^r)$ such that
$v_p(\sigma-1)=\alpha$).
\end{teo}
\begin{dimo}
First some remarks:

(1) All the functions in the above table (which are clearly class
functions) are distinct. In fact for functions belonging to different
rows, this follows from the fact that they have different degrees; for functions
of the first rows it's obvious; finally for the functions
$\psi^r_k\otimes\chi^r_k$ notice that
$$\psi^r_k\otimes\chi^r_k([z^{p^k}(1+p^k)])=p^{k-1}(p-1)\psi^r_k(1+p^k)$$
and so the difference follows from having chosen the $\psi^r_k$
among a rapresentative system of $G(p^r)^*/G(p^k)^*=<1+p^k>^*$ in
$G(p^r)^*$.\\
Hence, being all distinct, the number of these functions is
$$\#\{\mbox{Characters on the table}\}=p^{r-1}(p-1)+\sum_{k=1}^r
p^{r-k}=p^{r-1}(p-1)+\frac{p^r-1}{p-1}$$
which, for the Proposition 3.6, is equal to the number of conjugacy
classes of $K(p^r)$. So it's enough to show that
they are indeed irreducible characters of $K(p^r)$.

(2) The functions $\psi^r$ (hence also $\psi^r_k$) are irreducible
characters as they are induced by irreducible characters of $G(p^r)$
through the projection $K(p^r)\twoheadrightarrow G(p^r)$.

(3) As the $\psi^r$ (and so $\psi^r_k$) have values in
$\mathbb{S}^1=\{z\in \C\,:\,\mid z\mid=1\}$, in order to verify that
$\psi^r_k\otimes\chi^r_k$ are irreducible characters, it's enough to
verify that the $\chi^r_k$ are. In fact, by the remark (i) of the
theorem, the tensorial product of two irreducible characters is again  a
character; as for the irreducibility we can calculate their norm
(in the usual scalar product between characters) as follows
$$(\psi^r_k\otimes\chi^r_k,\psi^r_k\otimes\chi^r_k)_{K(p^r)}=\frac{1}
{|K(p^r)|}\sum_{g\in K(p^r)}\psi^r_k(g)\chi^r_k(g)\overline{\psi^r_k(g)}
\overline{\chi^r_k(g)}=$$
$$=\frac{1}{|K(p^r)|}\sum_{g\in K(p^r)}\chi^r_k(g)
\overline{\chi^r_k(g)}=(\chi^r_k,\chi^r_k)_{K(p^r)}$$
from which it follows that $\psi^r_k\otimes\chi^r_k$ is irreducible if and
only if $\chi^r_k$ is irreducible.

So these remarks tell us that to prove the theorem it's enough to show
that the functions $\chi^r_k$ are irreducible characters of $K(p^r)$.
We will do this first for the ``top'' function $\chi^r_r$ and then for
the others functions $\chi^r_k$, $1\leq k\leq r-1$, we will proceed by
induction on $r$.

$\boxed{\chi^r_r}$ \\
We will show first that it is a character and then that it is irreducible.

\underline{ CHARACTER }\\
Consider the primitive character $\chi$ of $C(p^r)$ which sends
$[1]_{C(p^r)}$ in $\z_{p^r}$ and induct with
respect to the inclusion $C(p^r)\hookrightarrow C(p^r)\rtimes G(p^r)=K(p^r)$.
The formula for the induced character $\chi^*$
(see \cite[Chapter XVIII, Section 6]{lang}) tells us
\begin{equation*}
\chi^*(z^i\sigma)=\sum_{\substack{\tau\in G(p^r)\\\tau z^i\sigma\tau^{-1}\in
    C(p^r)}}\chi(\tau z^i\sigma\tau^{-1})=
  \left\{\begin{aligned}
    &0\qquad&&\text{if }\sigma\neq 1,\\
    &\sum_{\tau\in G(p^r)}\chi(z^{i\tau})=\sum_{\tau\in
    G(p^r)}(\z_{p^r}^i)^{\tau}\:&&\text{if }\sigma=1.\\
  \end{aligned}\right.
\end{equation*}
Let us calculate the last summation.
\begin{lem}
Given $0\leq s\leq r$, we have
\begin{eqnarray}
\sum_{\tau\in G(p^r)}\z_{p^s}^{\tau}=
\left\{\begin{aligned}
&0&&\text{ if }&&s\geq2\\
&-p^{r-1}&&\text{ if }&&s=1\\
&p^{r-1}(p-1)\qquad&&\text{ if }&&s=0.\\
\end{aligned}\right.
\end{eqnarray}
\end{lem}
\begin{dimo}
For $t\geq s$ we have
\begin{eqnarray}
\sum_{x\in C(p^t)}\z_{p^s}^x=
\left\{\begin{aligned}
&0&&\text{ if }&&s\geq 1\\
&p^t\quad&&\text{ if }&&s=0\\
\end{aligned}\right. 
\end{eqnarray}
where the last one is obvious since $\z_{p^0}=1$ while the
first equality follows from the fact that $\z_{p^s}$,
with $1\leq s$, is a root of the polynomial $x^{p^s-1}+\cdots+1$
and so
$$\sum_{x\in C(p^t)}\z_{p^s}^x=\sum_{\substack{y\in C(p^s)\\z\in
C(p^{t-s})}} \z_{p^s}^{zp^s+y}=p^{t-s}\sum_{y\in C(p^s)}\z_{p^s}^y=0.$$
With the help of formula (3.5), we can write
\begin{equation*}
\sum_{\tau\in G(p^r)}\z_{p^s}^{\tau}=\sum_{x\in C(p^r)}\z_{p^s}^x-\sum_
{y\in C(p^{r-1})}\z_{p^s}^{py}=
\left\{\begin{aligned}
&0&&\text{ if }&&s\geq2\\
&-p^{r-1}&&\text{ if }&&s=1\\
&p^r-p^{r-1}=p^{r-1}(p-1)
\:&&\text{ if }&&s=0.\\
\end{aligned}\right.
\end{equation*}
\end{dimo}
Hence for $\chi^*$ we obtain
$$\chi^*([z^{p^{\beta}}\sigma_{\alpha}])=
\left\{\begin{aligned}
&0&&\text{ if }&&\sigma_{\alpha}\neq 1&&\text{ or }&&r-\beta\geq1\\
&-p^{r-1}&&\text{ if }&&\sigma_{\alpha}=1&&\text{ and }&&r-\beta=1\\
&p^{r-1}(p-1)&&\text{ if }&&\sigma_{\alpha}=1&&\text{ and }&&r-\beta=0\\
\end{aligned}\right.$$
which is precisely the definition of $\chi^r_r$. So, being induced
from a character of $C(p^r)$, $\chi^r_r$ is a character of $K(p^r)$.

\underline{IRREDUCIBILITY}\\
Now we calculate the scalar product of $\chi^r_r$ with itself. Since
$[z^{p^r}\sigma_r]$ contains only the identity and
$[z^{p^{r-1}}\sigma_r]$ contains $p-1$ elements, we have
$$(\chi^r_r,\chi^r_r)_{K(p^r)}=\frac{1}{|K(p^r)|}\sum_{g\in K(p^r)}\chi^r_r(g)
\overline{\chi^r_r(g)}=$$
$$=\frac{1}{|K(p^r)|}\left\{\left[p^{r-1}(p-1)\right]^2+
(p-1)\left[-p^{r-1}\right]^2\right\}=\frac{p^{2(r-1)}(p-1)p}{|K(p^r)|}=1$$
from which the irreducibility of $\chi^r_r$.

$\boxed{\chi^r_k,\, 1\leq k\leq r-1}$\\
Proceed by induction on $r$ (for $r=1$ we have only the function
$\chi^1_1$ which is an irreducible character for what proved before).
So let us assume, by induction hypothesis, that $\chi^{r-1}_k$,
$1\leq k\leq r-1$, are irreducible characters of  $K(p^{r-1})$
and let us show that $\chi^r_k$ is an irreducible character of $K(p^r)$.
In order to do this, consider the projection
$$\pi_r:K(p^r)\twoheadrightarrow K(p^{r-1})$$
obtained by reducing both $C(p^r)$ and  $G(p^r)$ modulo $p^{r-1}$.
Pull back the character $\chi^{r-1}_k$ to an irreducible character
$(\chi^{r-1}_k)'$ on $K(p^r)$. We will show
\begin{eqnarray}
(\chi^{r-1}_k)'=\chi^r_k 
\end{eqnarray}
and this will conclude the proof. Since, by definition,
 $$(\chi^{r-1}_k)'([z^{p^{\beta}}\sigma_{\alpha}])=
\chi^{r-1}_k(\pi_r([z^{p^{\beta}}\sigma_{\alpha}]))$$
and on the other hand
$$\pi_r([z^{p^{\beta}}\sigma_{\alpha}])=
\left\{\begin{aligned}
&\text{if}&&\alpha\leq r-1\:&&=[z^{p^{\beta}}\sigma_{\alpha}]\\
&\text{if}&&\alpha=r&&\left\{\begin{aligned}
                     &\text{and if}&&\beta\leq r-1\:&&=[z^{p^{\beta}}\cdot1]\\
                     &\text{and if}&&\beta=r\:&&=[1]\\
                             \end{aligned}\right.\\
\end{aligned}\right.$$
then we have
$$(\chi^{r-1}_k)'([z^{p^{\beta}}\sigma_{\alpha}])=
\left\{\begin{aligned}
&\text{if}&&\alpha\leq r-1&&\begin{sis}
          &\text{and if}&&\alpha<k&&\text{or}&&\beta<k-1&&=0\\
          &\text{and if}&&k\leq\alpha&&\text{and}&&\beta=k-1&&=-p^{k-1}\\
          &\text{and if}&&k\leq\alpha&&\text{and}&&k-1<\beta&&=p^{k-1}(p-1)\\
                             \end{sis}\\
&\text{if}&&\alpha=r&&\begin{sis}
                        &\text{and if}&&\beta<k-1&&=0\\
                        &\text{and if}&&\beta=k-1&&=-p^{k-1}\\
                        &\text{and if}&&k-1<\beta\leq r-1&&=p^{k-1}(p-1)\\
                        &\text{and if}&&\beta=r&&=p^{k-1}(p-1)\\
                       \end{sis}\\
\end{aligned}\right.$$
from which it follows that $(\chi^{r-1}_k)'=\chi^r_k$.
\end{dimo}
In the next section we will consider also the group $C(p^s)\rtimes G(p^r)$
for some $0\leq s\leq r$ where the semi-direct product is made with
respect to the map $G(p^r)\twoheadrightarrow
G(p^s) \cong \mbox{Aut}(C(p^s))$.
As a corollary of the preceding theorem, we now derive also explicit
formulas for the characters of this group (we can suppose $1\leq s\leq r$
because if $s=0$ we obtain the group $G(p^r)$ of which already we know
the characters). The notation used will be similar to that of the Theorem 3.7.
\begin{cor}
The irreducible characters of $C(p^s)\rtimes G(p^r)$ are
\begin{center}
\begin{tabular}{||p{8.5cm}||c|c||}
\hline
\hline
\bfseries CHARACTERS       & Number       & Degree\\
\hline
\hline
$\psi^r\qquad\text{ with }\psi^r\in G(p^r)^*$       &$p^{r-1}(p-1)$&$1$\\
\hline
$\psi^r_s\otimes\chi_s^s\;\text{with }\psi^r_s\mbox{ syst. of repr. for }
G(p^r)^*/G(p^s)^*$
& $p^{r-s}$ & $p^{s-1}(p-1)$\\
\hline
$\cdots$                  &$\cdots$      &$\cdots$ \\
\hline
$\psi^r_k\otimes\chi^s_k\;\text{with }\psi^r_k\mbox{ syst. of repr. for }
G(p^r)^*/G(p^k)^*$
& $p^{r-k}$ & $p^{k-1}(p-1)$\\
\hline
$\cdots$                 &$\cdots$      &$\cdots$ \\
\hline
$\psi^r_1\otimes\chi^s_1\;\text{with }\psi^r_1\mbox{ syst. of repr. for }
G(p^r)^*/G(p^1)^*$
& $p^{r-1}$ &$(p-1)$\\
\hline
\hline
\end{tabular}
\end{center}
where $\chi^s_k$, $1\leq k\leq s$, is the character so defined
$$\chi^s_k([z^{p^{\beta}}\sigma_{\alpha}])=
\left\{\begin{aligned}
&0\quad&&\text{ if }&&\alpha<k&&\text{ or }&&\beta<k-1,\\
&-p^{k-1}&&\text{ if }&&k\leq\alpha&&\text{ and }&&\beta=k-1,\\
&p^{k-1}(p-1)&&\text{ if }&&k\leq\alpha&&\text{ and }&&k-1<\beta.\\
\end{aligned}\right.$$
\end{cor}
\begin{dimo}
Observe first of all that the number of conjugacy classes of
$C(p^s)\rtimes G(p^r)$ is determined by the same rules of Theorem 2.4
except for the new condition $\beta=v_p(i)\leq s$. Hence the number of the
conjugacy classes can be calculated in this way
$$\#\{\mbox{Conjugacy classes}\}=$$
$$=\#\{[z\sigma]\}+\#\{[z^p\sigma]\,:\,[z^p\sigma]
\neq[z\sigma]\}+\cdots+\#\{[z^{p^s}\sigma]\,:\,[z^{p^s}\sigma]\neq[z^{p^{s-1}}
\sigma]\}=$$
$$=|G(p^r)|+|G(p^r)^1|+\cdots+|G(p^r)^s|=p^{r-1}(p-1)+p^{r-1}+\cdots+p^{r-s}=$$
\begin{equation}
=p^{r-1}(p-1)+p^{r-s}\frac{p^s-1}{p-1}.
\end{equation}
Now consider the projection
$$\pi:C(p^r)\rtimes G(p^r)\twoheadrightarrow C(p^s)\rtimes G(p^r)$$
obtained by reducing $C(p^r)$ modulo $p^s$.
From it, we deduce that the irreducible representations of
$C(p^s)\rtimes G(p^r)$ are exactly the representations of
$C(p^r)\rtimes G(p^r)$ which are the identity on $\ker\pi=<z^{p^s}>$.
This implies that an irreducible character
of $C(p^s)\rtimes G(p^r)$ induces, by composition with the projection $\pi$,
an irreducible character $\chi$ of $C(p^r)\rtimes G(p^r)$ such that
$$\chi_{|<z^{p^s}>}=\chi(1)\cdot 1_{|<z^{p^s}>}.$$
The only characters of $C(p^r)\rtimes G(p^r)$ which satisfy this property
are (with the notation of the Theorem 3.7):

(i) $\psi^r$, for which $\psi^r([z^{p^s}])=1=\chi([1])$;

(ii) $\psi^r_k\otimes \chi^r_k$, with  $1\leq k\leq s$, for which
$$\psi^r_k\otimes \chi^r_k([z^{p^s}]1)=p^{k-1}(p-1)=\psi^r_k\otimes
\chi^r_k([1]).$$
As their number is
$$p^{r-1}(p-1)+p^{r-1}+\cdots+p^{r-s}=p^{r-1}(p-1)+p^{r-s}\frac{p^s-1}{p-1}$$
which, for what observed at the beginning, is the number of the conjugacy
classes of $C(p^s)\rtimes G(p^r)$, necessarily they are all the
irreducible characters $C(p^s)\rtimes G(p^r)$. This proves the
theorem after having renamed $\chi^r_k$ as $\chi^s_k$.
\end{dimo}
In the next sections, in order to calculate the Artin conductor,
 we will be interested in knowing if the
restriction of a character to certain subgroups is trivial or not.
So we end this section with a result in this direction. First some
definitions.
\begin{defi}[Level]
We call level of a character of $C(p^s)\rtimes G(p^r)$ the number
so determined:
$$lev(\chi)=
\begin{sis}
&0&&\quad\text{if }\chi=\psi^r\in G(p^r)^*\\
&k&&\quad\text{if }\chi=\psi^r_k\otimes \chi^s_k.\\
\end{sis}$$
\end{defi}
\begin{defi}[Primitive degree]
It is called primitive degree (and indicated $pr$) of
$\psi^r\in G(p^r)^*$ the smallest number $0\leq\rho\leq r$ such that
$\psi^r$ is induced by a character of $G(p^{\rho})$ through the
projection
$G(p^r)\twoheadrightarrow G(p^{\rho}).$\\
It is called primitive degree of $\psi^r_k\otimes \chi^s_k$
the smallest number $k\leq\rho\leq r$ such that
$\psi^r_k$ is induced by a character of $G(p^{\rho})$ through the
projection
$G(p^r)\twoheadrightarrow G(p^{\rho}).$
\end{defi}
\begin{defi}[Null subgroup]
The null subgroup of a character $\chi$ of $C(p^s)\rtimes G(p^r)$
(indicated with $Gr(\chi)$) is the smallest subgroup of
$C(p^s)\rtimes G(p^r)$ such that
$$\chi_{| Gr(\chi)}=\chi(1)1_{| Gr(\chi)}$$
i.e. the corresponding representation is the identity.
\end{defi}
\begin{teo}
The null subgroup of a character $\chi$ of $C(p^s)\rtimes G(p^r)$
is equal to
$$ Gr(\chi)=C(p^{s-lev(\chi)})\rtimes G(p^r)^{pr(\chi)}.$$
\end{teo}
\begin{dimo}
Observe that:

(1) $\psi\in G(p^r)^*$ is equal to $1$ on $G(p^r)^t$ if and only if $t\geq
pr(\psi)$.

(2) $\chi^s_k([z^{p^{\beta}}\sigma])=\chi^s_k(1)=p^{k-1}(p-1)$ if and only if
$\beta\geq k$ and $v_p(\sigma-1)\geq k$.\\
From these two remarks it follows that
$$\chi([z^{p^{\beta}}\sigma])=\chi(1)\Rightarrow
\begin{sis}
&\beta\geq lev(\chi)\\
&\sigma\in G(p^r)^{pr(\chi)}\\
\end{sis}$$
and hence the theorem.
\end{dimo}

\section{Reduction to the prime power case}

In the section we begin to study the ramification of a prime $p$
in the extension $\Q(\z_m,\sqrt[m]{a})/\Q$ with the hypothesis $m$ odd and
if $p\mid m$ then $p\nmid v_p(a) \text{ or } p^{v_p(m)}\mid v_p(a)$.
The aim of this section is to show that the
wild part of the ramification is concentrated on the subextension
$\Q(\z_{p^r},\sqrt[p^r]{a})$, where $r=v_p(m)$, so that the higher ramification
groups can be calculated considering only this subextension.\\
First we want to determine which primes ramify in
$\Q(\z_m,\sqrt[m]{a})/\Q$.
\begin{lem}
Let $K$ be a number field and $L=K(\sqrt[n]{a})$.
If a prime $\p$ of $K$ doesn't divide $na$ then it is not ramified in $L$.
\end{lem}
\begin{dimo}
Consider the discriminant $d(L/K)$ of $L/K$. It holds:
$$d(L/K)\mid d_{L/K}(\sqrt[n]{a})\big|\left(N_{L/K}((x^n-a)'_{|x=\sqrt[n]{a}}
\right)=n^na^{n-1}.$$
hence if $\p\nmid (na)$, then $\p\nmid d(L/K)$ and so it is not
ramified.
\end{dimo}
\begin{cor}
The primes ramified in the extension $\Q(\z_m,\sqrt[m]{a})/\Q$ are
the divisor of $m$ or $a$.
\end{cor}
Next we study the primes which divides $a$ but not $m$.
\begin{teo}
If $p\nmid m$ then $\Q(\z_m,\sqrt[m]{a})/\Q$ has ramification
index respect to $p$ equal to
$$e(\Q(\z_m,\sqrt[m]{a})/\Q)=\frac{m}{(m,v_p(a))}.$$
In particular it is tamely ramified.
\end{teo}
\begin{dimo}
Consider the tower of extensions
$$\Q\subset \Q(\sqrt[m]{a})\subset \Q(\z_m,\sqrt[m]{a}).$$
The last extension is obtained by adding $\sqrt[m]{1}$ to the preceding one so that, as $p\nmid m$,
Lemma 4.1 implies that it is not ramified respect to $p$. So the ramification index
of the total extension is equal to the ramification index of the extension
$\Q\subset\Q(\sqrt[m]{a})$.\\
If we put $d=(v_p(a),m)$, then we can write
$$\begin{aligned}
&a&&=p^{d\alpha}a'\\
&m&&=dm'\\
\end{aligned}$$
with $(\alpha,m')=1$ and $p\nmid (dm'a')$. Now consider the tower of extensions
$$\Q\subset\Q(\sqrt[d]{a})\subset\Q(\sqrt[m]{a}).$$
Since $\sqrt[d]{a}=p^{\alpha}\sqrt[d]{a'}$, again by Lemma 4.1 we deduce that
$\Q\subset\Q(\sqrt[d]{a})=\Q(\sqrt[d]{a'})$ is not ramified respect to $p$.
Hence the total ramification index is equal to the ramification index of the extension
\begin{equation*}
\Q\left(\sqrt[d]{a'}\right)\subset\Q\left(\sqrt[m']{p^{\alpha}\sqrt[d]{a'}}
\right).
\end{equation*}
Since $(\alpha,m')=1$, there exist $s, t\in\Z$ such that $s\alpha-tm'=1$
with $(s,m')=1$. Now according to Theorem 2.3 we can transform the
extension as
$$\Q\left(\sqrt[m']{p^{\alpha}\sqrt[d]{a'}}\right)=\Q\left(\left(\sqrt[m']
{p^{\alpha}\sqrt[d]{a'}}\right)^sp^{-t}\right)=\Q\left(\sqrt[m']{p\sqrt[d]
{(a')^s}}\right)$$
and so, called $u=\sqrt[d]{(a')^s}$, we can complete with respect to the valuation $\p$-adic
($\p$ is one the primes of $\Q\left(\sqrt[d]{a'}\right)$ lying above $p$)
reducing ourselves to determine the ramification index of the local extension
\begin{equation*}
K\subset K(\sqrt[m']{pu})=L
\end{equation*}
where $K$ is a non ramified finite extension of $\Q_p$ and $u$ is an invertible element of $K$.
Look now to the valuation $v_L$-adic of the element $\sqrt[m']{pu}$:
$$v_L\left(\sqrt[m']{pu}\right)=\frac{v_L(pu)}{m'}=\frac{e(L/K)v_K(pu)}{m'}=
\frac{e(L/K)}{m'}$$
from which it follows that $m'\mid e(L/K)$; but since $e(L/K)\leq [L:K]\leq m'$
then $e(L/K)=m'$. The theorem follows from the definition of $m'=m/d=m/(m,v_p(a))$.
\end{dimo}
Now we come to the general case in which $p$ divides $m$ (and
possibly also $a$). If we write $m=p^rn$ with $p\nmid n$, then we
can split our extension as\\

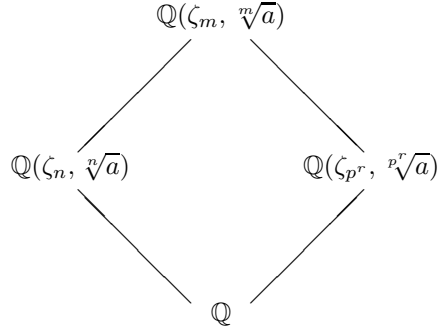
\setlength{\unitlength}{1mm}
\begin{picture}(80,40)(-20,0)
\linethickness{.5pt}
\put(40,40){\makebox(0,0)[t]{$\Q(\z_m,\sqrt[m]{a})$}}
\put(44,35){\line(1,-1){15}}
\put(60,20){\makebox(0,0)[t]{$\Q(\z_{p^r},\sqrt[p^r]{a})$}}
\put(36,35){\line(-1,-1){15}}
\put(20,20){\makebox(0,0)[t]{$\Q(\z_n,\sqrt[n]{a})$}}
\put(44,0){\line(1,1){15}}
\put(36,0){\line(-1,1){15}}
\put(40,0){\makebox(0,0)[t]{$\Q$}}
\end{picture}\\
\vspace{0.3mm}

Now we show how the determination of the ramification groups can be reduced to the study of the
extension $\Q(\z_{p^r},\sqrt[p^r]{a})$ that will be done in next sections.
\begin{teo}
The ramification index of $p$ in $\Q(\z_m,\sqrt[m]{a})/\Q$ ($m=p^rn$, with $p\nmid n$)
is equal to the following least common multiple
$$e(\Q(\z_m,\sqrt[m]{a})/\Q)=\left[\frac{n}{(n,v_p(a))},
e(\Q(\z_{p^r},\sqrt[p^r]{a})/\Q)\right]$$
while for the higher ramification groups we have
$$G(\Q(\z_m,\sqrt[m]{a})/\Q)^u=G(\Q(\z_{p^r},\sqrt[p^r]{a})/\Q)^u\:\:\text{ for } \: u>0$$
\end{teo}
\begin{dimo}
For the second assertion, observe that, since
$\Q(\z_n,\sqrt[n]{a})/\Q$ is tamely ramified, its ramification groups vanish for degree $>0$.
But this implies
$$G(\Q(\z_m,\sqrt[m]{a})/\Q)^u K(p^r)/K(p^r)=\{1\}\Rightarrow
G(\Q(\z_m,\sqrt[m]{a})/\Q)^u\subset K(p^r)$$
for $u>0$ and hence we conclude by taking the quotient with respect to $K(n)$.\\
For the first assertion, observe that question is local (so that we can take
the completion of all the fields involved respect to primes lying above $p$)
and that the preceding theorem tells us that the extension $\Q(\z_n,\sqrt[n]{a})$
is tamely ramified with index of ramification equal to $n/(n,v_p(a))$. So
the theorem descends from the following proposition.
\begin{pro}
Let $L_1/K$ and $L_2/K$ two disjoint finite extensions of $p$-local field with
ramification index $e_1$ and $e_2$ respectively.
If $L_1/K$ is tamely ramified (that is $p\nmid e_1$) then the
ramification $e$ of the composition $L_1L_2$ is equal to the lest common multiple of $e_1$ and $e_2$: 
$e=[e_1,e_2].$
\end{pro}
\begin{dimo}
Consider the maximal unramified subextensions $M_1$ and $M_2$ of respectively
$L_1$ and $L_2$ and look at the following diagram\\

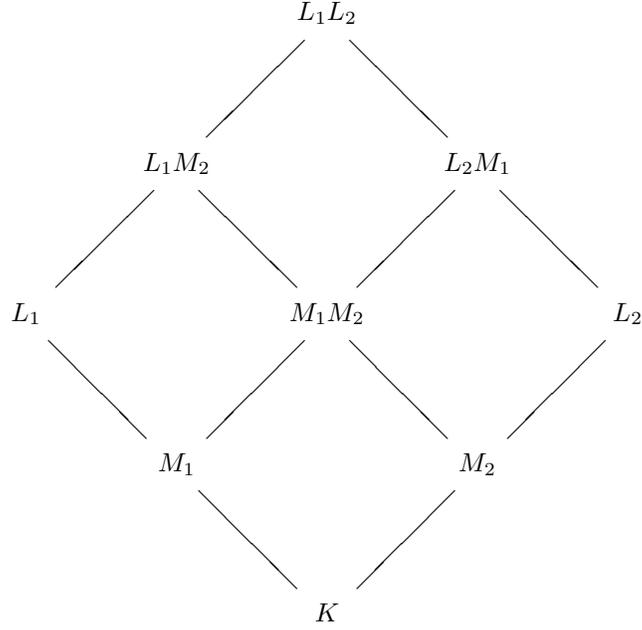
\setlength{\unitlength}{1mm}
\begin{picture}(0,80)(-20,0)
\linethickness{.5pt}
\put(40,80){\makebox(0,0)[t]{$L_1L_2$}}
\put(43,75){\line(1,-1){13}}
\put(37,75){\line(-1,-1){13}}
\put(60,60){\makebox(0,0)[t]{$L_2M_1$}}
\put(20,60){\makebox(0,0)[t]{$L_1M_2$}}
\put(63,55){\line(1,-1){13}}
\put(57,55){\line(-1,-1){13}}
\put(23,55){\line(1,-1){13}}
\put(17,55){\line(-1,-1){13}}
\put(0,40){\makebox(0,0)[t]{$L_1$}}
\put(40,40){\makebox(0,0)[t]{$M_1M_2$}}
\put(80,40){\makebox(0,0)[t]{$L_2$}}
\put(43,35){\line(1,-1){13}}
\put(77,35){\line(-1,-1){13}}
\put(3,35){\line(1,-1){13}}
\put(37,35){\line(-1,-1){13}}
\put(20,20){\makebox(0,0)[t]{$M_1$}}
\put(60,20){\makebox(0,0)[t]{$M_2$}}
\put(23,15){\line(1,-1){13}}
\put(57,15){\line(-1,-1){13}}
\put(40,0){\makebox(0,0)[t]{$K$}}
\end{picture}\\
\vspace{0.2mm}\\
From the property of stability of the unramified extensions
(\cite[Chapter 1, section 7]{cas}), it follows that $M_1M_2/K$ is unramified
while $L_1M_2/M_1M_2$ and $L_2M_1/M_2M_1$ are totally ramified
with
\begin{align*}
&e(L_1M_2/M_1M_2)=e(L_1/K),\\
&e(L_2M_1/M_2M_1)=e(L_2/K).
\end{align*}
So we can reduce precisely to the situation of the following lemma
and that will conclude the proof.
\begin{lem}
Let $M_1/M$ and $M_2/M$ two disjoint finite extensions of $p$-local field totally and
tamely ramified of degree $e_1$ and $e_2$. Then the composition $M_1M_2$ has
ramification index over $M$ equal to $[e_1,e_2]$.
\end{lem}
\begin{dimo}
According to the structure theorem for tamely totally ramified
extensions of local fields (see \cite[Chapter 1, section 8]{cas}), there exist in
$M$ elements $c_1$ e $c_2$ of valuation $1$ such that
\begin{align*}
&M_1=M(\sqrt[e_1]{c_1}),\\
&M_2=M(\sqrt[e_2]{c_2}).
\end{align*}
If we put $d=(e_1,e_2)$ then we can write
\begin{align*}
&e_1=de_1',\\
&e_2=de_2'.
\end{align*}
Consider the following diagram\\

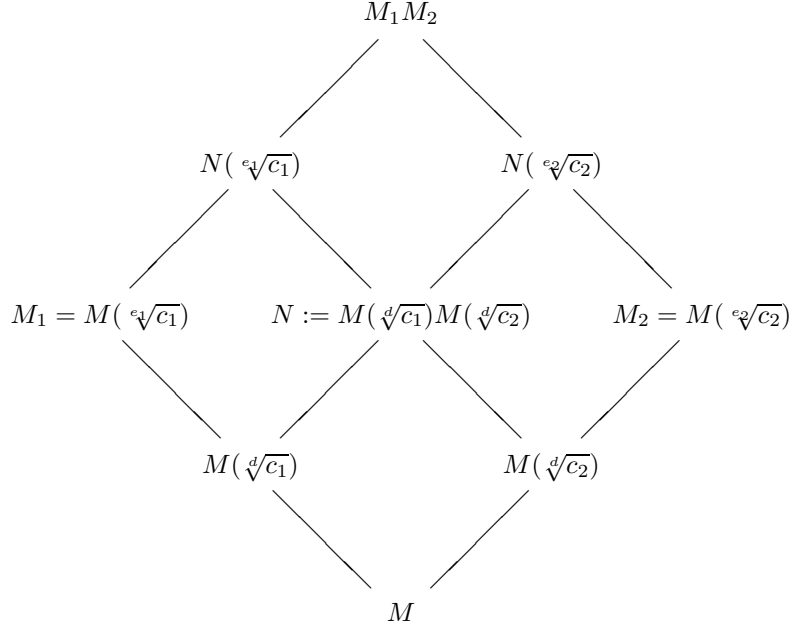
\setlength{\unitlength}{1mm}
\begin{picture}(0,80)(-20,0)
\linethickness{.5pt}
\put(40,80){\makebox(0,0)[t]{$M_1M_2$}}
\put(43,75){\line(1,-1){13}}
\put(37,75){\line(-1,-1){13}}
\put(60,60){\makebox(0,0)[t]{$N(\sqrt[e_2]{c_2})$}}
\put(20,60){\makebox(0,0)[t]{$N(\sqrt[e_1]{c_1})$}}
\put(63,55){\line(1,-1){13}}
\put(57,55){\line(-1,-1){13}}
\put(23,55){\line(1,-1){13}}
\put(17,55){\line(-1,-1){13}}
\put(0,40){\makebox(0,0)[t]{$M_1=M(\sqrt[e_1]{c_1})$}}
\put(40,40){\makebox(0,0)[t]{$N:=M(\sqrt[d]{c_1})M(\sqrt[d]{c_2})$}}
\put(80,40){\makebox(0,0)[t]{$M_2=M(\sqrt[e_2]{c_2})$}}
\put(43,35){\line(1,-1){13}}
\put(77,35){\line(-1,-1){13}}
\put(3,35){\line(1,-1){13}}
\put(37,35){\line(-1,-1){13}}
\put(20,20){\makebox(0,0)[t]{$M(\sqrt[d]{c_1})$}}
\put(60,20){\makebox(0,0)[t]{$M(\sqrt[d]{c_2})$}}
\put(23,15){\line(1,-1){13}}
\put(57,15){\line(-1,-1){13}}
\put(40,0){\makebox(0,0)[t]{$M$}}
\end{picture}\\
\vspace{0.2mm}\\
We will show that $e(N/M)=d$ and $e(M_1M_2/N)=e_1'e_2'$, from which it will follow that
$e(M_1M_2/M)=de_1'e_2'=[e_1,e_2]$ as requested.

\underline{ $e(N/M)=d$ }\\
Write $N$ as
$$N=M(\sqrt[d]{c_1},\sqrt[d]{c_2})=M(\sqrt[d]{c_1})\left(\sqrt[d]{\frac{c_1}
{c_2}}\right)$$
and since $v_M(c_1/c_2)=0$ Lemma 4.1 implies that the extension
$M(\sqrt[d]{c_1})\subset
M(\sqrt[d]{c_1})\left(\sqrt[d]{\frac{c_1}{c_2}}\right)$
is not ramified; hence $e(N/M)=e(M(\sqrt[d]{c_1})/M)=d$, q.e.d.

\underline{ $e(M_1M_2/N)=e_1'e_2'$ }\\
Observe that
$$\left\{\begin{aligned}
&M_1/M(\sqrt[d]{c_1}) &&\text{ tot. ram. of deg. }e_1'\\
&N/M(\sqrt[d]{c_1}) &&\text{ non ram. }\\
\end{aligned}\right.
\Rightarrow N(\sqrt[e_1]{c_1})/N \text{ tot. ram. of deg. }e_1'.$$
Analogously $N(\sqrt[e_2]{c_2})/N$ is totally ramified of degree $e_2'$.
Then, as $(e_1',e_2')=1$, the extension $M_1M_2= N(\sqrt[e_1]{c_1})
N(\sqrt[e_2]{c_2})$ is totally ramified over $N$ of degree $e_1'e_2'$, q.e.d. .
\end{dimo}
\end{dimo}
\end{dimo}
To summarize, we have shown in this section that we can
concentrate on the study of the ramification of $p$ in the
extension $\Q(\z_{p^r},\sqrt[p^r]{a})$. Besides, the original
hypothesis that $p^r\mid v_p(a) \text{ or } p\nmid v_p(a)$ splits into
the conditions $p\nmid a$ or $p\mid\mid a$.\\
In fact in the first case, we write $a=p^{p^s\alpha}a'$ with $p\nmid \alpha a'$ and
$s\geq r$ and obtain
$$\sqrt[p^r]{a}=p^{p^{s-r}\alpha}\sqrt[p^r]{a'}$$
so that replacing $a$ with $a'$ we can assume $p\nmid a$.\\
In the second case we can write $a=p^{v_p(a)}a'$ with $p\nmid
v_p(a)$. Hence there exist $s,t \in \Z$ such that $p^rs+v_p(a)t=1$ and $p\nmid t$,
which gives
$$a^t=\frac{pa'}{p^{p^rs}}$$
and by Theorem 2.3 we can replace $a$ by $pa'$.\\
So in the next two sections we will study the extension
$\Q(\z_{p^r},\sqrt[p^r]{a})$ distinguishing between these two cases.

\section{$\Q(\z_{p^r},\sqrt[p^r]{a})/\Q$ when $p\nmid a$}

First of all we want to complete our extensions. We recall the
following lemma of Kummer.
\begin{lem}[Kummer]
Let $L/K$ be an extension of number fields with rings of integer respectively $R_K$ and $R_L$.
Let $\theta\in R_L$ such that $L=K(\theta)$ and let $f(X)\in R_K[X]$ the minimal polynomial
of $\theta$ over $K$. Let $K_{\p}$ the completion of $K$ respect
to an ideal $\p$ of $R_K$ and let $R_{\p}$ its ring of integers.
If $f(X)$ factors on $R_{\p}[X]$ as
$$f(X)=\prod_{1\leq i\leq g}g_i(X)$$
then over $\p$ there are $g$ ideals of $R_L$ and the completions
with respect to this ideals are $K_{\p}(\theta_i)$
with $\theta_i$ root of $g_i(X)$.
\end{lem}
\begin{dimo}
See \cite[Chapter 2, section 10]{cas}.
\end{dimo}
If we apply this lemma to our situation we find the following result.
\begin{teo}
Let $0\leq s\leq r$ such that $a\in \Q_p^{p^{r-s}}$ and, if $s\neq 0$,
$a\not\in \Q_p^{p^{r-s+1}}$. In the extension $\Q(\z_{p^r},\sqrt[p^r]{a})/\Q$,
above $p$ there are $p^{r-s}$ prime ideals.
Besides, if $b\in \Q_p$ is such that $b^{p^{r-s}}=a$, then the
completion with respect to one of this ideals above $p$ is
$\Q_p(\z_{p^r},\sqrt[p^s]{b})$.
\end{teo}
\begin{dimo}
Consider the tower
$$\Q\subset\Q(\z_{p^r})\subset\Q(\z_{p^r},\sqrt[p^r]{a}).$$
We know that the cyclotomic extension
$\Q\subset \Q(\z_{p^r})$ is totally ramified at $p$.
Hence $\Q_p\subset\Q_p(\z_{p^r})$ is the completion with respect to $p$
and to the unique ideal over $p$.\\
Now apply the Kummer lemma to the extension
$\Q(\z_{p^r})(\sqrt[p^r]{a})$.
From our hypothesis on $s$, we know that $a=b^{p^{r-s}}$ with $b\in\Q_p$
such that, if $s\neq 0$, $b\not\in\Q_p^p$. The polynomial $X^{p^r}-a$
factorizes in $\Q_p(\z_{p^r})$ as
\begin{equation*}
X^{p^r}-a=(X^{p^s})^{p^{r-s}}-b^{p^{r-s}}=\prod_{i=1}^{p^{r-s}}
(X^{p^s}-\z_{p^{r-s}}^ib).\tag{*}
\end{equation*}
Now if $s=0$ then the polynomial splits into linear factors and so the
number of ideals over $p$ is $g=p^r$. Instead if $s>0$ then we know that
$b\not\in\Q_p^p$ and the Theorem 2.4 of Schinzel implies:
$$b\not\in\Q_p^p\Rightarrow  b\not\in\Q_p(\z_{p^r})^p$$
and so $\z_{p^{r-s}}^ib\not\in\Q_p(\z_{p^r})^p$, for every
$i=1,\cdots,p^{r-s}$. But then the Theorem 2.1 gives us that
$X^{p^s}-\z_{p^{r-s}}^ib$ is irreducible over $\Q_p(\z_{p^r})$ and so the (*)
is an irreducible factorization on $\Q_p(\z_{p^r})$. Now the lemma of Kummer
concludes the proof.
\end{dimo}
Now we want to compute the number $s$ of the preceding theorem. In
order to achieve this, we need the following theorem on the structure of
the units of $\Z_p$.
\begin{lem}
The group $U(\Z_p)$ of the invertible elements of $\Z_p$
admits the decomposition
$$U(\Z_p)\cong (\Z/p\Z)^*\times <1+p>.$$
Moreover the $p^k$-powers are
$$U(\Z_p)^{p^k}\cong (\Z/p\Z)^*\times <1+p^{k+1}>.$$
\end{lem}
\begin{dimo}
See \cite[p. 246-247]{has}.
\end{dimo}
Now the desired result.
\begin{teo}
Let $a\in \Q_p$ with $p\nmid a$. Then
$$a\in\Q_p^{p^r}\Leftrightarrow  p^{r+1}\mid \left(a^{p-1}-1\right).$$
Instead for $0<s\leq r$
$$a\in\Q_p^{p^{r-s}}\setminus \Q_p^{p^{r-s+1}}\Leftrightarrow
 p^{r-s+1}\mid\mid \left(a^{p-1}-1\right).$$
\end{teo}
\begin{dimo}
It follows from the preceding lemma that
$$a\in \Q_p^{p^k}\Leftrightarrow p^{k+1}\mid\left(a^{p-1}-1\right)$$
so that the conclusion is straightforward.
\end{dimo}
So we can restrict ourselves to the study of the local extension
$\Q_p(\z_{p^r},\sqrt[p^s]{b})$ with $0<s\leq r$, $p\nmid b$ and
$b\in\Q_p\setminus\Q_p^p$ (in fact the case $s=0$ reduces to the
cyclotomic extension $\Q_p(\z_{p^r})/\Q_p$ of which all is known).
Note that the fact that $b\not\in\Q_p^p$ is equivalent to say that
$p\mid\mid (b^{p-1}-1)$.

Now we want to determine the ramification degree of our extension.
\begin{teo}
The local extension $\Q_p(\z_{p^r},\sqrt[p^s]{b})/\Q_p$ with
$b\in\Q_p\setminus\Q_p^p$ is totally ramified.
\end{teo}
\begin{dimo}
Since it is well known that $\Q_p(\z_{p^r})/\Q_p$ is totally ramified,
it's enough to show that $\Q_p(\z_{p^r},\sqrt[p^s]{b})/\Q_p(\z_{p^r})$
is totally ramified, too. Suppose, on the contrary, that the inertia
degree $f$ is greater than $1$ and let
$$[\Q_p(\z_{p^r},\sqrt[p^s]{b})^{ur}:\Q_p(\z_{p^r})]=f$$
where $\Q_p(\z_{p^r},\sqrt[p^s]{b})^{ur}$ is the maximal unramified subextension of
$\Q_p(\z_{p^r},\sqrt[p^s]{b})$ over $\Q_p(\z_{p^r})$.
Then the unique subfield property (Theorem 2.2) gives
$$\Q_p(\z_{p^r},\sqrt[p^s]{b})^{ur}=\Q_p(\z_{p^r},\sqrt[f]{b})$$
and in particular $\sqrt[p]{b}\in \Q_p(\z_{p^r},\sqrt[p^s]{b})^{ur}$.
Now the structure theorem for the unramified extensions of local fields
(see \cite[Chapter 1, section 7]{cas}) says that $\Q_p(\z_{p^r})^{ur}$
is obtained from $\Q_p(\z_{p^r})$ adding
some roots of unity and so $\Q_p(\z_{p^r})^{ur}$is abelian over $\Q_p$.
In particular, since $\sqrt[p]{b}\in
\Q_p(\z_{p^r},\sqrt[p^s]{b})^{ur}$, $x^p-b$ has abelian Galois group
over $\Q_p$. Then the Theorem 2.4 of Schinzel gives $b=\gamma^p$
for some $\gamma\in \Q_p$ contrary to our hypothesis $b\not\in\Q_p^p$.
\end{dimo}
So our extension $\Q_p(\z_{p^r},\sqrt[p^s]{b})/\Q_p$ is totally
ramified of degree $\phi(p^r)p^s$ and, with the same proof of Theorem 3.2,
one shows that its Galois group is  $G=C(p^s)\rtimes G(p^r)$
where the semidirect product is made with respect to the map
$$G(p^r)\twoheadrightarrow
G(p^s) \cong \mbox{Aut}(C(p^s)).$$
Now we come to the heart of our work: the determination of the
higher ramification groups. The fundamental trick is to consider the
following filtration of subfields
$$\begin{array}{*{11}{c}}
\Q_p          & <  & \Q_p(\sqrt[p]{b})        & <  &\cdots&<   &\Q_p(\sqrt[p^s]{b})         &<   &\cdots&<   &\Q_p(\sqrt[p^r]{b})              \\

\triangle   &    &\triangle               &    &      &    &\wedge                    &    &      &    &\wedge              \\

\Q_p(\z_p)    &\lhd&\Q_p(\z_p,\sqrt[p]{b})    & <  &\cdots&<    &\Q_p(\z_p,\sqrt[p^s]{b})   &<   &\cdots&<   &\Q_p(\z_p,\sqrt[p^r]{b})              \\

\triangle   &    &\triangle               &    &      &     &\wedge                   &    &      &    &\wedge              \\

\Q_p(\z_{p^2})&\lhd&\Q_p(\z_{p^2},\sqrt[p]{b})&\lhd&\cdots&<   &\Q_p(\z_{p^2},\sqrt[p^s]{b})&<   &\cdots&<   &\Q_p(\z_{p^2},\sqrt[p^r]{b})\\

\triangle   &    &\triangle               &    &      &     &\wedge                   &    &      &    &\wedge                 \\

\vdots      &    &\vdots                  &    &      &    &\vdots                    &    &      &    &\vdots              \\

\triangle   &    &\triangle               &    &      &    &\triangle                 &    &      &    &\triangle                \\

\Q_p(\z_{p^r})&\lhd&\Q_p(\z_{p^r},\sqrt[p]{b})&\lhd&\cdots&\lhd&\Q_p(\z_{p^r},\sqrt[p^s]{b})&\lhd&\cdots&\lhd&\Q_p(\z_{p^r},\sqrt[p^r]{b})\\

\end{array}$$
where $\lhd$ means that the subfield on the right is normal over $\Q_p$,
where $<$ means that it is not normal over $\Q_p$. But observe that in
the second line every 1-step is normal, in the third line every 2-step is
normal and so on; while along the column all the extensions are normal.\\
Our strategy will be as follows: we shall determine the ramification on every step of the second
line and this, together with the knowledge of the ramification groups of
the cyclotomic extension (the first column of the diagram), will allow us
to calculate the total ramification groups of
$\Q_p(\z_{p^r},\sqrt[p^s]{b})/\Q_p$.
\begin{teo}
Consider the infinite tower of fields (in the hypothesis $p\mid\mid
(b^{p-1}-1)$)
$$K_0=\Q_p(\z_p)<K_1=\Q_p(\z_p,\sqrt[p]{b})<\cdots<K_i=\Q_p(\z_p,
\sqrt[p^i]{b})<\cdots$$
in which every extension is normal respect to the preceding one with
Galois group cyclic of order $p$. The unique break-number $b_i$ in the
ramification groups of $K_i/K_{i-1}$ is
$$b_i=1+p(p^{i-1}-1).$$
\end{teo}
\begin{dimo}
We shall first construct an uniformizer (i.e. an element of valuation
$1$) for $K_i/K_{i-1}$.
\begin{lem}
An uniformizer for the extension $K_i/K_{i-1}$ 
(with the hypothesis $p\mid\mid (b^{p-1}-1)$) is given by
\begin{equation*}
\pi_i=\frac{1-\z_p}{(b-\sqrt[p]{b})\cdots(\sqrt[p^{i-1}]{b}-\sqrt[p^i]{b})}.
\end{equation*}
\end{lem}
\begin{dimo}
Consider in $K_i$ the element $(\sqrt[p^{i-1}]{b}-\sqrt[p^i]{b})$.
For it we have
$$
N_{K_i/K_{i-1}}(\sqrt[p^{i-1}]{b}-\sqrt[p^i]{b})=\prod_{j=1}^{p}
(\sqrt[p^{i-1}]{b}-\z_p^j\sqrt[p^i]{b})=(\sqrt[p^{i-2}]{b}-\sqrt[p^{i-1}]{b}).
$$
So inductively we have that
\begin{equation*}
N_{K_i/K_0}(\sqrt[p^{i-1}]{b}-\sqrt[p^i]{b})=b^p-b
\end{equation*}
and this implies
\begin{equation*}
v_{K_i}(\sqrt[p^{i-1}]{b}-\sqrt[p^i]{b})=\frac{v_{K_0}(b^p-b)}{f(K_i/K_0)}=
v_{K_0}(b^p-b)=(p-1)v_p(b^p-b)=p-1 \tag{*}
\end{equation*}
in view of the hypothesis $p\mid\mid (b^{p-1}-1)$.
But this, together with the fact that $1-\z_p$ is a uniformizer
for $K_0=\Q_p(\z_p)/\Q_p$, allows to prove the theorem by induction
on $i\geq 1$. In fact:\\
\underline{ $i=1$ }\\
$$v_{K_1}\left(\frac{1-\z_p}{b-\sqrt[p]{b}}\right)=pv_{K_0}(1-\z_p)-v_{K_1}
(b-\sqrt[p]{b})=p-(p-1)=1.$$
\underline{ $i\Rightarrow i+1$ }\\
From (*) and the inductive hypothesis, we have
$$v_{K_i}\left(\frac{\pi_{i-1}}{\sqrt[p^{i-1}]{b}-\sqrt[p^i]{b}}\right)=p\cdot
v_{K_{i-1}}(\pi_{i-1})-v_{K_i}(\sqrt[p^{i-1}]{b}-\sqrt[p^i]{b})=p-(p-1)=1.$$
\end{dimo}
Now that we have an uniformizer $\pi_i$, to calculate the break-number of
$K_i/K_{i-1}$ we can simply let act on it the generator of the
group $\mbox{Gal}(K_i/K_{i-1})$\\
$\cong C(p)$ which sends
$\sqrt[p^i]{b}$ in $\z_p\sqrt[p^i]{b}$ (see \cite[Chapter IV, section 1, lemma 1]{ser}).
We obtain
$$1+b_i=v_{K_i}(s(\pi_i)-\pi_i)=v_{K_i}\left(\frac{\pi_{i-1}}
{\sqrt[p^{i-1}]{b}-\z_p\sqrt[p^i]{b}}-\frac{\pi_{i-1}}{\sqrt[p^{i-1}]{b}-
\sqrt[p^i]{b}}\right)=$$
$$=v_{K_i}(\pi_{i-1})+v_{K_i}\left(\frac{\sqrt[p^i]{b}\cdot(\z_p-1)}
{(\sqrt[p^{i-1}]{b}-\z_p\sqrt[p^i]{b})(\sqrt[p^{i-1}]{b}-\sqrt[p^i]{b})}
\right)=$$
$$=p\cdot v_{K_{i-1}}(\pi_{i-1})+v_{K_i}(\sqrt[p^i]{b})+p^i\cdot v_{K_0}
(\z_p-1)-2\cdot v_{K_i}(\sqrt[p^{i-1}]{b}-\sqrt[p^i]{b})=$$
$$=p+v_{K_i}(\sqrt[p^i]{b})+p^i-2(p-1)=p+p^i-2(p-1)$$
where in the last equality we have used that $p\nmid b$.
\end{dimo}
Now we can turn to the determination of the ramification groups of
our extension $\Q_p(\z_{p^r},\sqrt[p^s]{b})/\Q_p$.
\begin{nota}
In the next theorem we shall use the following notation

(i) for $i,s\leq r$, $C(p^s)\rtimes G(p^r)^i$ will indicate the
semidirect product made respect to the map
 $G(p^r)^i\hookrightarrow G(p^r)\twoheadrightarrow
G(p^s) \cong \mbox{Aut}(C(p^s))$.

(ii) with the numbers $d$ we shall indicate the difference of the
break-numbers in the ramification groups. Precisely $d_i$ will be
the difference between the $(i+1)$-th inferior break-number
and the $i$-th inferior break-number. Analogous meaning for $d^i$ respect
to the superior ramification groups.
\end{nota}
\begin{teo}
The ramification groups of $\Q_p(\z_{p^r},\sqrt[p^s]{b})/\Q_p$
($p\mid\mid (b^{p-1}-1)$) are:
$$\left\{\begin{aligned}
&G^0 &&= C(p^s)&&\rtimes&& G(p^r) && = G_0 \\
&G^{(i-1)+\frac{1}{p-1}}&&=C(p^{s-i+1})&&\rtimes&&
      G(p^r)^i&&=G_{\frac{2p^{2i-1}-p+1}{p+1}}\\
&G^i&&=C(p^{s-i})&&\rtimes&& G(p^r)^i&&=G_{\frac{(p-1)(p^{2i}-1)}{p+1}}  \\
&G^{s+j}&&=&&&&G(p^r)^{s+j}&&=G_{\frac{(p-1)(p^{2s}-1)}{p+1}+p^{2s}(p^j-1)}\\
\end{aligned}\right. .$$
with $1\leq i\leq s$ and $1\leq j\leq r-1-s$.
The difference of the break-numbers are:
$$\left\{\begin{aligned}
&d^0&&=0\\
&d^{2i-1}&&=1/(p-1)\\
&d^{2i}&&=(p-2)/(p-1)\\
&d^{2s+j}&&=1\\
\end{aligned}\right.
\quad\text{and}\quad
\left\{\begin{aligned}
&d_0&&=0\\
&d_{2i-1}&&=p^{2(i-1)}\\
&d_{2i}&&=p^{2i-1}(p-2)\\
&d_{2s+j}&&=p^{2s-1+j}(p-1).\\
\end{aligned}\right. $$
\end{teo}
The proof is by induction on $s$. For $s=0$ we reduce to a
cyclotomic extension and so the theorem follows from the
known ramification groups of the cyclotomic extensions of $\Q_p$.
(see \cite[Chapter IV, section 4]{ser}).
So we assume by induction that the theorem is true for $s<r$
and we will prove it for $s+1$.
\begin{lem}
The ramification groups of
$\Q_p(\z_{p^r},\sqrt[p^s]{b})/\Q_p(\z_p,
\sqrt[p^s]{b})$ are (call $\widetilde{G}=G(p^r)^1$
its Galois group):
$$\left\{\begin{aligned}
&\widetilde{G}^{(p-1)(p^i-1)}&&=G(p^r)^i&&=\widetilde{G}_{\frac{(p-1)
(p^{2i}-1)}{p+1}}\\
&\widetilde{G}^{(p-1)[(j+1)p^s-1]}&&=G(p^r)^{j+s}&&=\widetilde{G}_{\frac{(p-1)
(p^{2s-1}-1)}{p+1}+p^{2s}(p^j-1)}\\
\end{aligned}\right.$$
with $1\leq i\leq s$ and $1\leq j\leq r-1-s$. The differences of the break-numbers are
$$\left\{\begin{aligned}
&\widetilde{d}^{i-1}&&=p^{i-1}(p-1)^2\\
&\widetilde{d}^{s+j-1}&&=p^s(p-1)\\
\end{aligned}\right.
\quad\text{and}\quad
\left\{\begin{aligned}
&\widetilde{d}_{i-1}&&=p^{2i-2}(p-1)^2\\
&\widetilde{d}_{s+j-1}&&=p^{2s+j-1}(p-1).\\
\end{aligned}\right. $$
\end{lem}
\begin{dimo}
We know, by induction hypothesis, the ramification groups of\\
$\Q_p(\z_{p^r},\sqrt[p^s]{b})/\Q_p$
and so, to find the inferior ramification groups of\\
$\Q_p(\z_{p^r},\sqrt[p^s]{b})/\Q_p(\z_p,\sqrt[p^s]{b})$,
it's enough to intersect with $G(p^r)^1$. 
In particular for the inferior $d$ we obtain
$$\left\{\begin{aligned}
&\widetilde{d}_{i-1}&&=d_{2i-1}+d_{2i}=p^{2(i-1)}+p^{2i-1}(p-2)=p^{2i-2}
                       (p-1)^2\\
&\widetilde{d}_{s+j-1}&&=d_{2s+j}=p^{2s+j-1}(p-1).\\
\end{aligned}\right.$$
Now we can pass to the superior $d$:
$$\left\{\begin{aligned}
&\widetilde{d}^{i-1}&&=\frac{\widetilde{d}_{i-1}}{p^{i-1}}&&=p^{i-1}(p-1)^2\\
&\widetilde{d}^{s+j-1}&&=\frac{\widetilde{d}_{s+j-1}}{p^{s+j-1}}&&=p^s(p-1).\\
\end{aligned}\right.$$
And now from the $d_i$ and $d^i$ it's easy to calculate the inferior
and superior ramification groups.
\end{dimo}
Now using Theorem 5.6 we can move to the field
$\Q_p(\z_{p^r},\sqrt[p^{s+1}]{b})$ leaving fixed the base
$\Q_p(\z_p,\sqrt[p^s]{b})$.
\begin{lem}
The ramification groups of $\Q_p(\z_{p^r},\sqrt[p^{s+1}]{b})/\Q_p(\z_p,
\sqrt[p^s]{b})$ are (call ${}^sG=C(p)\rtimes G(p^r)^1$ its Galois group):
$$\left\{\begin{aligned}
&{}^sG^{(p-1)(p^i-1)}&&=C(p)&&\rtimes&& G(p^r)^i&&={}^sG_{\frac{(p-1)(p^{2i}-1)}
{p+1}}\\
&{}^sG^{1+p(p^s-1)}&&=C(p)&&\rtimes &&G(p^r)^{s+1}&&={}^sG_{\frac{p^{2s+1}-p+1}
{p+1}}\\
&{}^sG^{(p-1)[(j+1)p^s-1]}&&=&& && G(p^r)^{s+j}&&={}^sG_{\frac{(p-1)(p^{2s+2}-1)}{p+1}+p^{2s+2}(p^{j-1}-1)}\\
\end{aligned}\right.$$
with $1\leq i\leq s$ and $1\leq j\leq r-s-1$.
The differences of the break-numbers are:
$$\left\{\begin{aligned}
&{}^sd^{i-1}&&=p^{i-1}(p-1)^2\\
&{}^sd^s&&=p^s\\
&{}^sd^{s+1}&&=p^s(p-2)\\
&{}^sd^{s+j}&&=p^s(p-1)\\
\end{aligned}\right.
\quad\text{and}\quad
\left\{\begin{aligned}
&{}^sd_{i-1}&&=p^{2i-2}(p-1)^2&\quad&\\
&{}^sd_s&&=p^{2s}&&\\
&{}^sd_{s+1}&&=p^{2s+1}(p-2)&&\\
&{}^sd_{s+j}&&=p^{2s+j}(p-1)&&\text{ for }2\leq j.\\
\end{aligned}\right. $$
\end{lem}
\begin{dimo}
Consider the diagram
$$\begin{array}{cccl}
\Q_p(\z_p,\sqrt[p^s]{b})&\lhd&\Q_p(\z_p,\sqrt[p^{s+1}]{b})& \\
\triangle&&\triangle&G(p^r)^1\\
\Q_p(\z_{p^r},\sqrt[p^s]{b})&\lhd&\Q_p(\z_{p^r},\sqrt[p^{s+1}]{b})& \\
&C(p)&&\\
\end{array}$$
where the fixed fields of the subgroups $C(p)$ and $G(p^r)^1$
are respectively $\Q_p(\z_{p^r},\sqrt[p^s]{b})$
and $\Q_p(\z_p,\sqrt[p^{s+1}]{b})$.
From Theorem 5.6 and Lemma 5.9, we see that the unique break-numbers in
the quotients are:
$$\left\{\begin{aligned}
&{}^sG^{1+p(p^s-1)}/G(p^r)^1=C(p)\\
\end{aligned}\right.
\quad\text{and}\:
\left\{\begin{aligned}
&{}^sG^{(p-1)(p^i-1)}/C(p)&&=G(p^r)^i\\
&{}^sG^{(p-1)[(j+1)p^s-1]}/C(p)&&=G(p^r)^{s+j}.\\
\end{aligned}\right.$$
Since $(p-1)(p^s-1)<1+p(p^s-1)<(p-1)(2p^s-1)$, we conclude that

(1) The superior ramification groups of order $>1+p(p^s+1)$ are contained
in $G(p^r)^1$ and so ${}^sG^{(p-1)[(j+1)p^s-1]}=G(p^r)^{s+j}$.

(2) ${}^sG^{1+p(p^s-1)}=C(p)\rtimes G(p^r)^{s+1}$.

(3) For the ramification groups of order $<1+p(p^s+1)$ we have
    ${}^sG^{(p-1)(p^i-1)}$\\
    $=C(p)\rtimes G(p^r)^i$.\\
Next, with easily calculations (similar to that of the preceeding lemma), 
one can pass to superior $d$, inferior $d$ and finally to inferior 
ramification groups.
\end{dimo}
Now, always using Theorem 5.6, we can go down with the base field.
\begin{lem}
The ramification groups of $\Q_p(\z_{p^r},\sqrt[p^{s+1}]{b})/\Q_p(\z_p,
\sqrt[p^{s-t}]{b})$ with $0\leq t\leq s$ are (call
${}^{s-t}G=C(p^{t+1})\rtimes G(p^r)^1$):
$$\left\{\begin{aligned}
&{}^{s-t}G^{(p-1)(p^i-1)}&&=C(p^{t+1})&&\rtimes&&G(p^r)^i\\
&{}^{s-t}G^{p^{s-t+1}+(p-1)[kp^{s-t}-1]}&&=C(p^{t+1-k})&&\rtimes&&G(p^r)^
{s-t+1+k}\\
&{}^{s-t}G^{(p-1)[(k+2)p^{s-t}-1]}&&=C(p^{t-k})&&\rtimes&&G(p^r)^{s-t+1+k}\\
&{}^{s-t}G^{(p-1)[(r+j+1)p^{s-t}-1]}&&=&& &&G(p^r)^{s+j}\\
\end{aligned}\right.$$
with $1\leq i\leq s-t$, $0\leq k\leq t$ and $2\leq j\leq r-s-1$.
The differences of the break-numbers are:
$$\left\{\begin{aligned}
&{}^{s-t}d^{i-1}&&=p^{i-1}(p-1)^2\\
&{}^{s-t}d^{s-t+2k}&&=p^{s-t}\\
&{}^{s-t}d^{s-t+1+2k}&&=p^{s-t}(p-2)\\
&{}^{s-t}d^{s+t+j}&&=p^{s-t}(p-1)\\
\end{aligned}\right.
\quad\text{ and }\:
\left\{\begin{aligned}
&{}^{s-t}d_{i-1}&&=p^{2i-2}(p-1)^2\\
&{}^{s-t}d_{s-t+2k}&&=p^{2(s-t)+2k}\\
&{}^{s-t}d_{s-t+1+2k}&&=p^{2(s-t)+2k+1}(p-2)\\
&{}^{s-t}d_{s+t+j}&&=p^{2s+j}(p-1).\\
\end{aligned}\right.$$
\end{lem}
\begin{dimo}
The proof is by induction on $t$. For $t=0$ it reduces to Lemma 5.10.
So let us assume the theorem true for $t<s$ and prove it for $t+1$.\\
Consider the tower of subfields
$$\Q_p(\z_p,\sqrt[p^{s-t-1}]{b})\lhd\Q_p(\z_p,\sqrt[p^{s-t}]{b})\lhd
\Q_p(\z_{p^r},\sqrt[p^{s+1}]{b}).$$
Theorem 5.6 tells us that the break-number of
$\Q_p(\z_p,\sqrt[p^{s-t}]{b})/\Q_p(\z_p,\sqrt[p^{s-t-1}]{b})$ is
$$d_{s-t}=1+p(p^{s-t-1}-1).$$
and since $(p-1)(p^{s-t-1}-1)<1+p(p^{s-t-1}-1)<(p-1)(p^{s-t}-1)$,
we must distinguish for the differences of break-numbers the following
cases:

(1) The first $(s-t-1)$ break-number remains the same as the initial case
 $t$ and the corresponding ramification groups have an increase of the
first factor from $C(p^{t+1})$ to  $C(p^{t+2})$, i.e.
$$1\leq i\leq s-t\Rightarrow
\left\{\begin{aligned}
&{}^{s-t-1}d^{i-1}&&={}^{s-t}d^{i-1}=p^{i-1}(p-1)^2\\
&{}^{s-t-1}d_{i-1}&&={}^{s-t}d_{i-1}=p^{2i-2}(p-1)^2\\
&{}^{s-t-1}G^{(p-1)(p^i-1)}&&=C(p^{t+2})\rtimes G(p^r)^i.\\
\end{aligned}\right.$$

(2) It appears a new break-number for which one has
$$\left\{\begin{aligned}
&{}^{s-t-1}G^{1+p(p^{s-t-1}-1)}={}^{s-t-1}G^{p^{s-t}+(p-1)[0p^{s-t-1}-1]}=
C(p^{t+2})\rtimes G(p^r)^{s-t+1}\\
&{}^{s-t-1}d^{s-t-1}=1+p(p^{s-t-1}-1)-(p-1)(p^{s-t-1}-1)=p^{s-t-1}\\
&{}^{s-t-1}d_{s-t+1}=p^{2(s-t-1)}.\\
\end{aligned}\right.$$

(3) For the successive break-number it holds
$$\left\{\begin{aligned}
&{}^{s-t-1}d_{s-t+2}={}^{s-t}d_{s-t-1}-{}^{s-t-1}d_{s-t-1}=p^{2(s-t-1)+1}
(p-2)\\
&{}^{s-t-1}d^{s-t+2}=\frac{{}^{s-t-1}d_{s-t+2}}{p^{s-t+2}}=p^{s-t-1}(p-2)\\
&{}^{s-t-1}G^{2(p-1)p^{s-t-1}-(p-1)}=C(p^{t+1})\rtimes G(p^r)^{s-t-1}.\\
\end{aligned}\right.$$

(4) For the remaining break-number one has (for $1\leq k+1\leq t+1$ and
$2\leq j\leq r-1-s$):
$$\left\{\begin{aligned}
&{}^{s-t}d_{s-t+2k}&&={}^{s-t-1}d_{s-t+2k+1}&&={}^{s-t-1}d_{s-(t+1)+2(k+1)}\\
&{}^{s-t}d_{s-t+1+2k}&&={}^{s-t-1}d_{s-t+1+2k+1}&&={}^{s-t-1}d_{s-(t+1)+1+
2(k+1)}\\
&{}^{s-t}d_{s+j+t}&&={}^{s-t-1}d_{s+j+(t+1)}&&\\
\end{aligned}\right.$$
because of the appearance of the new ramification group. And for
the same reason one has that for every $h\geq 2(s-t)$
$${}^{s-t-1}d^{h+1}=\frac{{}^{s-t}d^h}{p}$$
from which we conclude that
$$\left\{\begin{aligned}
&{}^{s-t-1}d^{s-(t+1)+2(k+1)}&&=p^{s-(t+1)}\\
&{}^{s-t-1}d^{s-(t+1)+1+2(k+1)}&&=p^{s-(t+1)}(p-2)\\
&{}^{s-t-1}d^{s+(t+1)+j}&&=p^{s-(t-1)}(p-1).\\
\end{aligned}\right.$$
Now it's easy to compute the superior ramification
groups
$$\left\{\begin{aligned}
&{}^{s-t-1}G^{p^{s-t}+(p-1)[(k+1)p^{s-t}-1]}&&=C(p^{t+2-(k+1)})&&\rtimes&&G(p^r)
^{s-(t+1)+1+(k+1)}\\
&{}^{s-t-1}G^{(p-1)[(k+3)p^{s-t-1}-1]}&&=C(p^{t+1-(k+1)})&&\rtimes&&G(p^r)^
{s-(t+1)+1+(k+1)}\\
&{}^{s-t-1}G^{(p-1)[(t+1+j+1)p^{s-t-1}-1]}&&=&& &&G(p^r)^{s+j}.\\
\end{aligned}\right.$$
\end{dimo}
Putting $t=s$ in the preceding lemma, we obtain the following corollary
\begin{cor}
The ramification groups of $\Q_p(\z_{p^r},\sqrt[p^{s+1}]{b})/\Q_p(\z_p)$
are (call ${}^0G=C(p^{s+1})\rtimes G(p^r)^1$ its Galois group):
$$\left\{\begin{aligned}
&{}^0G^{1+k(p-1)}&&=C(p^{s+1-k})&&\rtimes&&G(p^r)^{k+1}\\
&{}^0G^{(k+1)(p-1)}&&=C(p^{s-k})&&\rtimes&&G(p^r)^{k+1}\\
&{}^0G^{(s+j)(p-1)}&&=&& &&G(p^r)^{s+j}\\
\end{aligned}\right.$$
where $0\leq k\leq s$ and $2\leq j\leq r-s-1$.
The differences of the break-numbers are
$$\left\{\begin{aligned}
&{}^0d^{2k}&&=1\\
&{}^od^{1+2k}&&=p-2\\
&{}^0d^{2s+j}&&=p-1\\
\end{aligned}\right.
\quad\text{ and }\quad
\left\{\begin{aligned}
&{}^0d_{2k}&&=p^{2k}\\
&{}^0d_{1+2k}&&=p^{1+2k}(p-2)\\
&{}^0d_{2s+j}&&=p^{2s+j}(p-1).\\
\end{aligned}\right.$$
\end{cor}
\begin{dimo}[of the Theorem 5.8]
Now we can conclude our proof by induction showing that the
theorem is true for $s+1$.
Consider the composition
$$\Q_p\lhd\Q_p(\z_p)\lhd\Q_p(\z_{p^r},\sqrt[p^{s+1}]{b}).$$
Since $\Q_p(\z_p)$ is the maximal subextensions tamely ramified,
the inferior and superior ramification groups
$\Q_p\lhd\Q_p(\z_{p^r},\sqrt[p^{s+1}]{b})$
differ from those of
$\Q_p(\z_p)\lhd\Q_p(\z_{p^r},\sqrt[p^{s+1}]{b})$
only for the $0$-th group which now is
$G_0=C(p^{s+1})\rtimes G(p^r)$. Hence for the inferior $d$ it holds
$$\left\{\begin{aligned}
&d_0&&=0\\
&d_{2k+1}&&={}^0d_{2k}=p^{2k}\\
&d_{2k+1}&&={}^0d_{2k+1}=p^{2k+1}(p-2)\\
&d_{2s+j+1}&&={}^0d_{2s+j}=p^{2s+j}(p-1)\\
\end{aligned}\right.$$
with $0\leq k\leq s$ and $2\leq j\leq r-s-1$.
Analogously for the superior $d$ one has
$$\left\{\begin{aligned}
&d^0&&=0\\
&d^{2k+1}&&=\frac{1}{p-1}\\
&d^{2k+1}&&=\frac{p-2}{p-1}\\
&d^{2s+j+1}&&=1.\\
\end{aligned}\right.$$
Now we calculate the superior ramification  groups:
$$\left\{\begin{aligned}
&G^0&&=C(p^{s+1})&&\rtimes&&G(p^r)\\
&G^{k+\frac{1}{p-1}}&&=C(p^{s+1-k})&&\rtimes&&G(p^r)^{k+1}\\
&G^{k+1}&&=C(p^{s-k})&&\rtimes&&G(p^r)^{k+1}\\
&G^{s+j}&&=&& &&G(p^r)^{s+j}\\
\end{aligned}\right.$$
which, after putting $i=k+1$ and $j'=j-1$, gives us the desired groups.
Now its easy to calculate the inferior ramification groups.
\end{dimo}
Now that we know the ramification groups of our extension
$\Q_p(\z_{p^r},\sqrt[p^s]{b})/\Q_p$, we can
calculate the local Artin conductor of the characters
of $C(p^s)\rtimes G(p^r)$. For this we shall use the formula
(\cite[Chapter VI, section 2, exercise 2]{ser})
$$f(\chi)=\chi(1)(1+c^{\chi})$$
where $c^{\chi}$ is the biggest real $u$ such that $\chi_{|G^u}\neq
\chi(1)\cdot id_{|G^u}$ (if $\chi=id$ we put $c^{id}=-1$).
\begin{teo}
For a character $\chi$ of $C(p^s)\rtimes G(p^r)$, we have
$$c^{\chi}=
\begin{sis}
&pr(\chi)-1&&\:\text{ if }lev(\chi)< pr(\chi)\:\text{ or }\:0=lev{\chi},\\
& pr(\chi)-1+\frac{1}{p-1}&&\:\text{ if }\: 0<lev(\chi)=pr(\chi),\\
\end{sis}$$
(here $lev(\chi)$ has the meaning of Def. 3.10).
\end{teo}
\begin{dimo}
From the definition 3.12 of the null subgroup of a character, we have that
$c^{\chi}$ (with $\chi\neq id$) is the biggest real number such that
$G^{c^{\chi}}\not\subset Gr(\chi).$
Now from Theorem 5.8 it follows that
$$G^u\not\subset Gr(\chi)=C(p^{s-lev(\chi)})\rtimes G(p^r)^
{pr(\chi)}$$
$$\Leftrightarrow
\begin{sis}
&u\leq lev(\chi)-1+\frac{1}{p-1}&&\:\text{ if }0<lev(\chi),\\
&\qquad\qquad\text{ or }\\
&u\leq pr(\chi)-1&&\:\text{ if }0<pr(\chi),\\
\end{sis}$$
and hence if $0<lev(\chi)$
$$c^{\chi}=\max\left\{lev(\chi)-1+\frac{1}{p-1},pr(\chi)-1\right\}=$$
$$=\begin{sis}
& pr(\chi)-1&&\:\text{ if }0<lev(\chi)<pr(\chi),\\
& pr(\chi)-1+\frac{1}{p-1}&&\:\text{ if } 0<lev(\chi)=pr(\chi).\\
\end{sis}$$
\end{dimo}
\begin{cor}
The local Artin conductor of $\chi$ is
$$f(\chi)=
\begin{sis}
&p^{lev(\chi)-1}(p-1) pr(\chi) &&\:\text{ if }
 0<lev(\chi)< pr(\chi),\\
&p^{lev(\chi)-1}(p-1)\left[ pr(\chi)+\frac{1}{p-1}\right]&&\:\text{ if }
 0<pr(\chi)=lev(\chi),\\
&pr(\chi)&&\:\text{ if } lev(\chi)=0.
\end{sis}$$
\end{cor}
\begin{dimo}
It follows at once from the preceding theorem and from the formula
$$f(\chi)=\chi(1)(1+c^{\chi}).$$
\end{dimo}
As an application of the formula for the local Artin conductor, we now
calculate the $p$-component of the discriminant of the extension
$\Q(\z_{p^r},\sqrt[p^r]{a})/\Q$. First of all, since above $p$ there are
$p^{r-s}$ prime ideals, we can complete our extension and apply the
formula
\begin{eqnarray}
v_p\left[d(\Q(\z_{p^r},\sqrt[p^r]{a})/\Q)\right]=p^{r-s}
v_p\left[d(\Q_p(\z_{p^r},\sqrt[p^s]{b})/\Q_p)\right]. 
\end{eqnarray}
For the discriminant of the extension
$\Q_p(\z_{p^r},\sqrt[p^s]{b})/\Q_p$, we use the local
conductor-discriminant formula (\cite[Chapter 6, section 3]{ser})
\begin{eqnarray}
v_p(d)=\sum_{\chi\in G^*}\chi(1)f(\chi). 
\end{eqnarray}
\begin{teo}
For the extension $\Q_p(\z_{p^r},\sqrt[p^s]{b})/\Q_p$ with $p\nmid b$ we
have
$$v_p\left[d(\Q_p(\z_{p^r},\sqrt[p^s]{b})/\Q_p)\right]=p^s\left[rp^r-(r+1)p^{r-1}
\right]+2\frac{p^{2s}-1}{p+1}.$$
\end{teo}
First of all we need the following lemma.
\begin{lem}
If $G=C(p^s)\rtimes G(p^r)$ then we have for its characters the two
formulas
$$
\text{ Card }\{\chi\in G^*:\:lev(\chi)=0 \:\text{ and }\:0\leq pr(\chi)=t\leq r\}=
\begin{sis}
&1&&\:\text{ if }\:t=0,\\
&p-2&&\:\text{ if }\: t=1,\\
&p^{t-2}(p-1)^2&&\:\text{ if } t\geq 2;\\
\end{sis}
$$
$$
\text{ Card }\{\chi\in G^*:\:lev(\chi)=k>0\:\text{ and }\:k\leq pr(\chi)=t\leq
r\}=
\begin{sis}
&1&&\:\text { if } \:t=k,\\
&(p-1)p^{t-k-1}&&\;\text{ if }\: t>k.\\
\end{sis}
$$
\end{lem}
\begin{dimo}
The first formula follows from:
$$
\text{ Card }\{\chi\in G^*:\:lev(\chi)=0 \:\text{ and }\:0\leq pr(\chi)=t\leq r\}=
|G(p^t)|-|G(p^{t-1})|
$$
while the second follows from (remember that we take only a system
of representatives for $G(p^r)/G(p^k)$):
$$
\text{ Card }\{\chi\in G^*:\:lev(\chi)=k>0\:\text{ and }\:k\leq pr(\chi)=t\leq
r\}=\frac{|G(p^t)|-|G(p^{t-1})|}{|G(p^k)|}.
$$
\end{dimo}
\begin{dimo}[of the Theorem 5.15]
According to Corollary 5.14 we can split the sum in the
formula (5.2) as
$$
v_p(d)=\sum_{lev(\chi)=0}pr(\chi)+\sum_{0<lev(\chi)}
p^{2(lev(\chi)-1)}(p-1)^2 pr(\chi)+
$$
\begin{eqnarray}
+\sum_{0<lev(\chi)=pr(\chi)}p^{2(lev(\chi)-1)}(p-1)^2
\left[\frac{1}{p-1}\right]
\end{eqnarray}
and with the help of the preceding lemma we can calculate these three 
sums.
Let us begin with the first sum:
\begin{eqnarray}
\sum_{lev(\chi)=0}pr(\chi)=(p-2)\cdot 1+\sum_{t=2}^r p^{t-2}(p-1)^2t=
rp^r-(r+1)p^{r-1}. 
\end{eqnarray}
For the second sum:
$$
\sum_{0<lev(\chi)}p^{2(lev(\chi)-1)}(p-1)^2pr(\chi)=
$$
$$
=\sum_{k=1}^s p^{2(k-1)}(p-1)^2k+
\sum_{k=1}^s\sum_{t=k+1}^rp^{2(k-1)}(p-1)^2t(p-1)p^{t-k-1}=
$$
\begin{eqnarray}
=(p^s-1)\left[rp^r-(r+1)p^{r-1}\right]+\frac{p^{2s}-1}{p+1}. 
\end{eqnarray}
Finally for the third sum we have:
\begin{eqnarray}
\sum_{0<lev(\chi)=pr(\chi)}p^{2(lev(\chi)-1)}(p-1)^2\left[\frac{1}{p-1}\right]=
\sum_{k=1}^sp^{2(k-1)}(p-1)=\frac{p^{2s}-1}{p-1}. 
\end{eqnarray}
Now adding the formulas (5.4), (5.5) and (5.6) we obtain what asserted.
\end{dimo}
So for the global extension we have:
\begin{cor}
The $p$-adic valuation of the discriminant of
$\Q(\z_{p^r},\sqrt[p^r]{a})/\Q$ with $p\nmid a$ is equal to
$$
p^r\left[rp^r-(r+1)p^{r-1}\right]+2\frac{p^{r+s}-p^{r-s}}{p+1}.
$$
\end{cor}

\section{$\Q_p(\z_{p^r},\sqrt[p^r]{a})/\Q_p$ in the case $p\mid\mid a$}

The other case of our interest is the extension $\Q_p(\z_{p^r},\sqrt[p^r]{a})/\Q_p$
in the case $p\mid\mid a$. First we want to complete this
extension.
\begin{lem}
In the extension $\Q(\z_{p^r},\sqrt[p^r]{a})/\Q$ with $p\mid\mid a$,
above $p$ there is only one prime ideal. So the completion of
our extension with respect to $p$ is
$\Q_p(\z_{p^r},\sqrt[p^r]{a})/\Q_p$.
\end{lem}
\begin{dimo}
It follows from Theorem 5.2 since if $p\mid\mid a$ then $a\not\in \Q_p^p$.
\end{dimo}
\begin{cor}
The local extension $\Q_p(\z_{p^r},\sqrt[p^r]{a})/\Q_p$ with $p\mid\mid a$
is totally ramified.
\end{cor}
\begin{dimo}
It follows from the preceding lemma together with Theorem 5.5.
\end{dimo}
Now we come to the determination of the higher ramification groups. We begin with an
analogue of Theorem 5.6.
\begin{teo}
Consider the tower of extensions (in the hypothesis $p\mid\mid a$)
$$K_0=\Q_p(\z_p)<K_1=\Q_p(\z_p,\sqrt[p]{a})<\cdots<K_i=\Q_p(\z_p,
\sqrt[p^i]{a})<\cdots$$
in which every extension is normal respect to the preceding one with
Galos group cyclic of order $p$. The unique break-number $b_i$ in the
ramification groups of $K_i/K_{i-1}$ is
$$b_i=p^i.$$
\end{teo}
\begin{dimo}
We shall first construct an uniformizer for $K_i/K_{i-1}$.
\begin{lem}
An uniformizer for the extension $K_i/K_{i-1}$
(with the hypothesis $p\mid\mid a$) is given by
\begin{equation*}
\pi_i=\frac{1-\z_p}{\sqrt[p]{a}\cdots\sqrt[p^i]{a}}.
\end{equation*}
\end{lem}
\begin{dimo}
Consider in $K_i$ the element $\sqrt[p^i]{a}$.
Thanks to the hypothesis $p\mid\mid a$
\begin{equation*}
v_{K_i}(\sqrt[p^i]{a})=\frac{v_{K_i}(a)}{p^i}=\frac{p^iv_{K_0}(a)}{p^i}=p-1.\tag{*}
\end{equation*}
But this, together with the fact that $1-\z_p$ is a uniformizer
for $K_0=\Q_p(\z_p)/\Q_p$, allows to prove the theorem by induction
on $i\geq 1$. In fact:\\
\underline{ $i=1$ }\\
$$v_{K_1}\left(\frac{1-\z_p}{\sqrt[p]{a}}\right)=pv_{K_0}(1-\z_p)-v_{K_1}
(\sqrt[p]{a})=p-(p-1)=1.$$
\underline{ $i\Rightarrow i+1$ }\\
From (*) and the inductive hypothesis, we have
$$v_{K_i}\left(\frac{\pi_{i-1}}{\sqrt[p^i]{a}}\right)=p\cdot
v_{K_{i-1}}(\pi_{i-1})-v_{K_i}(\sqrt[p^i]{a})=p-(p-1)=1.$$
\end{dimo}
Now that we have an uniformizer $\pi_i$, to calculate the break-number $b_i$ of
$K_i/K_{i-1}$ we can simply let act on it the generator of the
group $\mbox{Gal}(K_i/K_{i-1})\cong C(p)$ which sends
$\sqrt[p^i]{a}$ in $\z_p\sqrt[p^i]{a}$ (see \cite[Chapter IV, section 1, lemma 1]{ser}).
We obtain
$$1+b_i=v_{K_i}(s(\pi_i)-\pi_i)=v_{K_i}\left(\frac{\pi_{i-1}}
{\z_p\sqrt[p^i]{a}}-\frac{\pi_{i-1}}{\sqrt[p^i]{a}}\right)=$$
$$=v_{K_i}(\pi_{i})+v_{K_i}\left(\frac{1-\z_p}{\z_p}\right)
=1+p^iv_{K_0}(1-\z_p)=1+p^i.$$
\end{dimo}
Now that we have an analogue of the Theorem 5.6 we can calculate
the higher ramification groups of our extension with the same
inductive method that we used in the preceding section. So, even
if we are now interested only in the extension
$\Q_p(\z_{p^r},\sqrt[p^r]{a})/\Q_p$,
we must consider also the intermediate extensions
$\Q_p(\z_{p^r},\sqrt[p^s]{a})/\Q_p$ with $0\leq s\leq r$.
\begin{teo}
The ramification groups of $\Q_p(\z_{p^r},\sqrt[p^s]{a})/\Q_p$ are:
$$\left\{\begin{aligned}
&G^0&&=C(p^s)&&\rtimes&&G(p^r)\\
&G^1&&=C(p^s)&&\rtimes&&G(p^r)^1\\
&G^{i+\frac{1}{p-1}}&&=C(p^{s-i+1})&&\rtimes&&G(p^r)^{i+1}&&\\
&G^{i+1}&&=C(p^{s-i})&&\rtimes&&G(p^r)^{i+1}\\
&G^{s+j+1}&&=&&&&G(p^r)^{s+j+1}&&\\
\end{aligned}\right. $$
with $1\leq i\leq s$ and $1\leq j\leq r-2-s$ (with the convention that if
$r-2\leq s$ the last row doesn't appear).
\end{teo}
As the proof follows identically the same steps of that of
Theorem 5.8, we give only the statement of the various lemmas
involved without any proof (they are in fact conceptually
identical to those relative to Theorem 5.8).\\
By induction, we assume the theorem true for a certain $s<r$
and prove it for $s+1$ (the base of induction is as always the
cyclotomic extension).
\begin{lem}
The ramification groups of $\Q_p(\z_{p^r},\sqrt[p^s]{a})/\Q_p(\z_p,
\sqrt[p^s]{a})$ are (call $\widetilde{G}=G(p^r)^1$ its Galois group):
$$\left\{\begin{aligned}
&\widetilde{G}^{p-1}&&=G(p^r)^1\\
&\widetilde{G}^{(p-1)p^i}&&=G(p^r)^{i+1}\\
&\widetilde{G}^{(j+1)(p-1)p^s}&&=G(p^r)^{j+s+1}\\
\end{aligned}\right.$$
with $1\leq i\leq s$ and $1\leq j\leq r-2-s$.\\
\end{lem}
Now using Theorem 6.3 we can move to the field $\Q_p(\z_{p^r},
\sqrt[p^{s+1}]{a})$ leaving fixed the base $\Q_p(\z_p,\sqrt[p^s]{a})$.
\begin{lem}
The ramification groups of $\Q_p(\z_{p^r},\sqrt[p^{s+1}]{a})/\Q_p(\z_p,
\sqrt[p^s]{a})$ are (call ${}^sG=C(p)\rtimes G(p^r)^1$ its Galois group):
$$\left\{\begin{aligned}
&{}^sG^{p-1}&&=C(p)&&\rtimes&&G(p^r)^1\\
&{}^sG^{(p-1)p^i}&&=C(p)&&\rtimes&&G(p^r)^{i+1}\\
&{}^sG^{p^{s+1}}&&=C(p)&&\rtimes &&G(p^r)^{s+2}\\
&{}^sG^{(j+1)(p-1)p^s}&&=&& && G(p^r)^{s+j+1}\\
\end{aligned}\right.$$
with $1\leq i\leq s$ and $1\leq j\leq r-s-2$.
\end{lem}
Now again using Theorem 6.3 we can go down with the base field.
\begin{lem}
The ramification groups of $\Q_p(\z_{p^r},\sqrt[p^{s+1}]{a})/\Q_p(\z_p,
\sqrt[p^{s-t}]{a})$ with $0\leq t\leq s$ are (call ${}^{s-t}G=C(p^{t+1})\rtimes
G(p^r)^1$ its Galois group):
$$\left\{\begin{aligned}
&{}^{s-t}G^{p-1}&&=C(p^{t+1})&&\rtimes&&G(p^r)^1\\
&{}^{s-t}G^{(p-1)p^i}&&=C(p^{t+1})&&\rtimes&&G(p^r)^{i+1}\\
&{}^{s-t}G^{p^{s-t}\left[(2k+1)p-2k\right]}&&=C(p^{t+1-k})&&\rtimes&&G(p^r)^
{s-t+2+k}\\
&{}^{s-t}G^{(k+2)(p-1)p^{s-t}}&&=C(p^{t-k})&&\rtimes&&G(p^r)^{s-t+2+k}\\
&{}^{s-t}G^{(j+t+1)(p-1)p^{s-t}}&&=&& &&G(p^r)^{s+j+1}\\
\end{aligned}\right.$$
with $1\leq i\leq s-t$, $0\leq k\leq t$ and $2\leq j\leq r-s-2$.
\end{lem}
When $t=s$ the preceding lemma gives the following corollary.
\begin{cor}
The ramification groups of  $\Q_p(\z_{p^r},\sqrt[p^{s+1}]{a})/\Q_p(\z_p)$
are (call ${}^0G=C(p^{s+1})\rtimes G(p^r)^1$ its Galois group):
$$\left\{\begin{aligned}
&{}^0G^{p-1}&&=C(p^{s+1})&&\rtimes&&G(p^r)^1\\
&{}^0G^{p+2k(p-1)}&&=C(p^{s+1-k})&&\rtimes&&G(p^r)^{k+2}\\
&{}^0G^{(k+2)(p-1)}&&=C(p^{s-k})&&\rtimes&&G(p^r)^{k+2}\\
&{}^0G^{(s+j+1)(p-1)}&&=&& &&G(p^r)^{s+j+1}\\
\end{aligned}\right.$$
with $0\leq k\leq s$ and $2\leq j\leq r-s-2$.
\end{cor}
\begin{dimo}[of Theorem 6.5]
Now the theorem follows by observing that the inferior ramification groups of
$\Q_p(\z_{p^r},\sqrt[p^{s+1}]{a})/\Q_p$ differ from those of \\
$\Q_p(\z_{p^r},\sqrt[p^{s+1}]{a})/\Q_p(\z_p)$ only in degree $0$
where we must substitute $G(p^r)^1$ with $G(p^r)$.
\end{dimo}
We include in a corollary the case $s=r$ that interests us.
\begin{cor}
The ramification groups of $\Q_p(\z_{p^r},\sqrt[p^r]{a})/\Q_p$
are: $$\left\{\begin{aligned}
&G^0&&=C(p^r)&&\rtimes&&G(p^r)&&=G_0\\
&G^{1}&&=C(p^{r})&&\rtimes&&G(p^r)^1&&=G_{(p-1)}\\
&G^{i+\frac{1}{p-1}}&&=C(p^{r-i+1})&&\rtimes&&G(p^r)^{i+1}&&=G_{\frac{2p^{2i}+p-1}{p+1}}\\
&G^{i+1}&&=C(p^{r-i})&&\rtimes&&G(p^r)^{i+1}&&=G_{(p-1)\frac{p^{2i}+1}{p+1}}\\
&G^{r-1+\frac{1}{p-1}}&&=C(p^2)&& && &&=G_{\frac{2p^{2r-2}+p-1}{p+1}}\\
&G^{r+\frac{1}{p-1}}&&=C(p)&& && &&=G_{\frac{p^{2r}+p^{2r-2}+p-1}{p+1}}\\
\end{aligned}\right. $$
with $1\leq i\leq r-2$.
\end{cor}
Now that we have the higher ramification groups, we can
calculate the local Artin conductor of the characters of
$C(p^r)\rtimes G(p^r)$ using again the formula (\cite[Chapter VI, section 2, exercise 2]{ser})
$$f(\chi)=\chi(1)(1+c^{\chi})$$
where $c^{\chi}$ is the biggest ramification group such that
$$\chi_{|G^{c^{\chi}}}\neq \chi(1)1_{|G^{c^{\chi}}}.$$
\begin{teo}
For a charcter $\chi$ of $C(p^r)\rtimes G(p^r)$ we have
$$c^{\chi}=
\begin{sis}
& pr(\chi)-1&&\:\text{ if } lev(\chi)+2\leq pr(\chi)\:\text{ or }
 0=lev{\chi},\\
& lev(\chi)+\frac{1}{p-1}&&\:\text{ if }
 0<lev(\chi)\leq pr(\chi)\leq lev(\chi)+1.\\
\end{sis}$$
\end{teo}
\begin{dimo}
From the definition of the null subgroup of a character (see definition 4.12),
we deduce that $c^{\chi}$ is the biggest real number such that
$G^{c^{\chi}}\not\subset Gr(\chi).$
So from Corollary 6.10 we deduce
$$G^u\not\subset Gr(\chi)=C(p^{r-lev(\chi)})\rtimes G(p^r)^
{ pr(\chi)}$$
$$\Leftrightarrow
\begin{sis}
&u\leq k+\frac{1}{p-1}&&\:\text{ if  }0<lev(\chi)\\
&\qquad\qquad\text{ or }\\
&u\leq pr(\chi)-1&&\:\text{ if } 0<pr(\chi).\\
\end{sis}$$
Hence if $0<lev(\chi)$ we have
$$c^{\chi}=\max\left\{lev(\chi)+\frac{1}{p-1},
pr(\chi)-1\right\}=$$
$$=\begin{sis}
& pr(\chi)-1&&\:\text{ if } lev+2\leq pr(\chi),\\
& lev(\chi)+\frac{1}{p-1}&&\:\text{ if } pr(\chi)\leq lev(\chi)+1.\\
\end{sis}$$
\end{dimo}
\begin{cor}
The local Artin conductor of $\chi$ is
 $$f(\chi)=
\begin{sis}
&pr(\chi)&&\:\text{ if } lev(\chi)=0,\\
&p^{lev(\chi)-1}(p-1)pr(\chi)&&\:\text{ if  }2<lev(\chi)+2\leq pr(\chi),\\
&p^{lev(\chi)-1}(p-1)[lev(\chi)+1+\frac{1}{p-1}]&&\:\text{ if }
 0<lev(\chi)\leq pr(\chi)\leq lev(\chi)+1.\\
\end{sis}$$
\end{cor}
\begin{dimo}
It follows from the preceding theorem and the formula
$f(\chi)=\chi(1)(1+c^{\chi})$.
\end{dimo}
Let us now calculate the $p$-component of the discriminant of the extension
$\Q(\z_{p^r},\sqrt[p^r]{a})/\Q$ with $p\mid\mid a$. As it is a totally ramified extension
we have
$$v_p(d(\Q(\z_{p^r},\sqrt[p^r]{a})/\Q))=v_p(d(\Q_p(\z_{p^r},\sqrt[p^r]{a})/
\Q_p))$$
so that we can focus ourselves on the local extension
$\Q_p(\z_{p^r},\sqrt[p^r]{a})/\Q_p$.
\begin{teo}
If $p\mid\mid a$ then we have
$$v_p\left(d(\Q_p(\z_{p^r},\sqrt[p^r]{a})/\Q_p)\right)=
rp^{2r-1}(p-1)+p\frac{p^{2r}-1}{p+1}-p\frac{p^{2r-3}+1}{p+1}.$$
\end{teo}
\begin{dimo}
Using the local conductor-discriminant formula and Corollary 6.12
we have:
$$
v_p(d)=\sum_{\chi\in G^*}\chi(1)f(\chi)=\sum_{lev(\chi)=0}pr(\chi)+
\sum_{lev(\chi)>0}p^{2(lev(\chi)-1)}(p-1)^2pr(\chi)+
$$
$$
+\sum_{0<lev(\chi)=pr(\chi)-1}p^{2(lev(\chi)-1)}(p-1)^2\frac{1}{p-1}
+\sum_{0<lev(\chi)=pr(\chi)}p^{2(lev(\chi)-1)}(p-1)^2
\left[1+\frac{1}{p-1}\right].
$$
The first two sums are identical to those in the formulas (5.4) and (5.5) 
of the fourth section. For the third sum we have:
$$
\sum_{0<lev(\chi)=pr(\chi)-1}p^{2(lev(\chi)-1)}(p-1)^2\frac{1}{p-1}=
\sum_{k=1}^{r-1}p^{2(lev(\chi)-1)}(p-1)^2(p-1)\frac{1}{p-1}=
$$
\begin{eqnarray}
=(p-1)\frac{p^{2r-2}-1}{p+1}. 
\end{eqnarray}
Finally the fourth sum is equal to
$$\sum_{0<lev(\chi)=pr(\chi)}p^{2(lev(\chi)-1)}\left[1+\frac{1}{p-1}\right]=
\sum_{k=1}^rp^{2(k-1)}(p-1)^2\left(1+\frac{1}{p-1}\right)=
$$
\begin{eqnarray}
=p\frac{p^{2r}-1}{p+1}. 
\end{eqnarray}
Adding together the formulas (5.4), (5.5), (6.1) and (6.2) we have the 
desired result.
\end{dimo}

\end{document}